\documentclass[11pt]{article}
   \usepackage{amsfonts}
\usepackage{latexsym,tikz}
\usepackage{amsmath}
 
\title{Multiple Dedekind Zeta Functions}
\author{Ivan Emilov Horozov}
\date{February 2, 2008}

\newcommand \nc {\newcommand}
\nc \proof {\noindent {\em{Proof.\/ }}} \nc \qed {$\Box$\hfill}
\newtheorem{theorem}{Theorem}[section]
\newtheorem{lemma}[theorem]{Lemma}
\newtheorem{proposition}[theorem]{Proposition}
\newtheorem{corollary}[theorem]{Corollary}
\newtheorem{definition}[theorem]{Definition}
\newtheorem{example}[theorem]{Example}
\newtheorem{remark}[theorem]{Remark}
\newtheorem{conjecture}[theorem]{Conjecture}
\newtheorem{question}[theorem]{Question}
\nc \bth[1] {\begin{theorem}\label{t#1} } \nc \ble[1]
{\begin{lemma}\label{l#1} } \nc \bpr[1]
{\begin{proposition}\label{p#1} } \nc \bco[1]
{\begin{corollary}\label{c#1} } \nc \bde[1]
{\begin{definition}\label{d#1}\rm } \nc \bex[1]
{\begin{example}\label{e#1}\rm } \nc \bre[1]
{\begin{remark}\label{r#1}\rm } \nc \bcon[1]
{\begin{conjecture}\label{con#1}\rm } \nc \bque[1]
{\begin{question}\label{que#1}\rm }
\nc {\eth} { \end{theorem} } \nc {\ele} { \end{lemma} } \nc
{\epr}{ \end{proposition} } \nc {\eco} { \end{corollary} } \nc
{\ede} {\end{definition} } \nc {\eex} { \end{example} } \nc {\ere}
{\end{remark} } \nc {\econ} { \end{conjecture} } \nc {\eque}
{\end{question} }
\def\a{\alpha}

\def\e{\epsilon}

\def \di {{\mathrm d}}

\def \N {{\mathbb N}}
\def \Z {{\mathbb Z}}
\def \Q {{\mathbb Q}}
\def \R {{\mathbb R}}
\def \C {{\mathbb C}}

\def \Re { {\mathrm{Re}} }
\begin{document}

\title{{\LARGE\bf{Multiple Dedekind Zeta Functions}}}

\author{
Ivan ~Horozov
\thanks{E-mail: horozov@math.wustl.edu}
\\ \hfill\\ \normalsize \textit{Department of Mathematics,}\\
\normalsize \textit{Washington University in St. Louis,}\\
\normalsize \textit{One Brookings Drive, Campus box 1146}\\
\normalsize \textit {Saint Louis, MO 63130, USA }  \\
}
\date{}
\maketitle

\begin{abstract}
In this paper we define multiple Dedekind zeta values (MDZV),
using a new type of iterated integrals, called iterated integrals
on a membrane. One should consider MDZV as a number theoretic generalization of Euler's multiple zeta values. Over imaginary quadratic fields MDZV capture, in particular, multiple Eisenstein series \cite{ZGK}. We give an analogue of multiple Eisenstein series over real quadratic field and an alternative definition of values of multiple Eisenstein-Kronecker series \cite{G2}. Each of them is a special case of multiple Dedekind zeta values.
MDZV are interpolated into functions that we call multiple Dedekind zeta functions (MDZF). We show that MDZF have integral representation, can be written as infinite sum, and have analytic continuation. We compute explicitly the value of a multiple residue of certain MDZF over a quadratic number  field at the point $(1,1,1,1)$. Based on such computations, we state two conjectures about MDZV.
\end{abstract}

\tableofcontents
\setcounter{section}{-1}

\section{Introduction}
\label{Intro}
Multiple Dedekind zeta functions generalize Dedekind zeta
functions in the same way the multiple zeta functions generalize 
the Riemann zeta function. Let us recall known definitions of the
above functions. The Riemann zeta function is defined as
$$\zeta(s)=\sum_{n>0}\frac{1}{n^s},$$
where $n$ is an integer.
Multiple zeta functions are defined as
$$\zeta(s_1,\dots,s_m)=\sum_{0<n_1<\dots<n_m}\frac{1}{n_1^{s_1}\dots n_m^{s_m}},$$
where $n_1,\dots,n_m$ are integers.
Special values of the Riemann zeta function $\zeta(k)$ and of the multiple zeta functions 
$\zeta(k_1,\dots,k_m)$ were defined by Euler \cite{Eu}.
The Riemann zeta function is closely related to the ring of integers.

Dedekind zeta function $\zeta_K(s)$ is an analogue of the Riemann
zeta function, which is closely related to the algebraic integers
${\cal{O}}_K$ in a number field $K$. It is defined as
$$\zeta_K(s)=\sum_{\mathfrak{a}\neq(0)}\frac{1}{N_{K/\Q}(\mathfrak{a})^s},$$
where the sum is over all ideals $\mathfrak{a}$ different from the zero ideal $(0)$ and
$N(\mathfrak{a})=\#|{\cal{O}}_K/\mathfrak{a}|$ is the norm of the ideal $\mathfrak{a}$.

A definition of multiple Dedekind zeta functions should combine
ideas from multiple zeta functions and from Dedekind zeta
functions. 

There is a definition of multiple Dedekind zeta functions due to Masri \cite{Mas}.
Let us recall  his definition. Let $K_1,\dots,K_m$ be number fields and let
${\cal{O}}_{K_i}$, for $i=1,\dots,m$, be the corresponding rings of
integers.  Let $\mathfrak{a}_i$, for $i=1,\dots,m$, be ideals in
${\cal{O}}_{K_i}$, respectively. Then he defines
$$\zeta(K_1,\dots,K_d; s_1,\dots s_m)
=
\sum_{0<N(\mathfrak{a}_1)<\dots<N(\mathfrak{a}_m)}
\frac{1}{N(\mathfrak{a}_1)^{s_1}\dots N(\mathfrak{a}_m)^{s_m}}.$$

We propose a different definition. The advantage of our definition
is that it leads to more properties: analytic, topological and
algebraic-geometric. Let us give an explicit formula for a multiple
Dedekind zeta function, in a case when it is easier to formulate.
Let $K$ be a number field with ring of integers ${\cal{O}}_K$. Let
$U_K$ be the group of units in ${\cal{O}}_K$. Let $C$ be a cone
inside of a fundamental domain of ${\cal{O}}_K$ modulo $U_K$. (More
precisely, $C$ has to be a positive unimodular simple cone as defined in
Section \ref{sec cone}. A fundamental domain for ${\cal{O}}_K$
modulo $U_K$ can be written as a finite union of unimodular simple
cones.) For such a cone $C$, we define a multiple Dedekind zeta
function
\begin{align}
&\zeta_{K;C}(s_1,\dots,s_1;\dots;s_m,\dots,s_m)=\\
&\sum_{\alpha_1,\dots,\alpha_m\in C}
\frac{1}{N(\alpha_1)^{s_1}N(\alpha_1+\alpha_2)^{s_2}\cdots N(\alpha_1+\dots+\alpha_m)^{s_m}}.
\end{align}

The key new ingredient in the definition of multiple Dedekind zeta
functions is the definition of iterated integrals on a membrane.
This is a higher dimensional analogue of iterated path integrals.
In the iterated integrals on a membrane the iteration
happens in $n$-directions.
Such iterated integrals were defined in \cite{ModSym}
generalizing Manin's
non-commutative modular symbol \cite{Man} to higher dimensions in
some cases, essentially for Hilbert modular surfaces. 

{\bf{Structure of the paper:}}

In Subsection \ref{subsec examples polylog}, we recall definitions of multiple zeta values and of polylogarithms by giving many explicit formulas. In Subsection \ref{subsec examples polylog Gauss int}, we generalize the previous formulas to multiple Dedekind zeta values over the Gaussian integers via many examples.

In Section \ref{sec membrane}, we give two Definitions of iterated
integrals on a membrane. 
The first definition is more intuitive. It can be used to generalize the first few formulas for MDZV over the Gaussian integers from Subsection \ref{subsec examples polylog Gauss int}.
The second Definition is the one needed for 
the definition of multiple Dedekind zeta values. It is needed in order to express special values of the multiple Eisenstein series via MDZV, when the modular parameter has a value in an imaginary quadratic field.

In Section \ref{sec cone}, we use some basic algebraic number
theory (see \cite{IR}), in order to construct the functions that
we integrate. We use an idea of Shintani (see \cite{Sh}, \cite{C})
for defining a cone. We associate a product of
geometric series to every unimodular simple cone. This is the type
of functions that we integrate. Lemma \ref{lemma fundamental
domain}  shows that a fundamental domain for the non-zero integer
${\cal{O}}_K-\{0\}$  modulo the units $U_K$ can be written as a
finite union of unimodular simple cones.

In Section \ref{sec MDZF}, we define Dedekind polylogarithms
associated to a positive unimodular simple cone. Theorem \ref{thm DZ}
expresses Dedekind zeta values in terms of Dedekind
polylogarithms. The heart of the section is Definition \ref{def
MDZV} of multiple Dedekind zeta values (MDZV) as an iterated integral over a membrane and Definition \ref{def MDZF} of multiple Dedekind zeta functions (MDZF) in terms of an integral representation. 
Theorems \ref{thm MDZV} and \ref{thm MDZF} express MDZV and MDZF as an infinite sum.
At the end of the Section \ref{sec MDZF}, we give many examples.  
Examples 1 and 2 are the simplest multiple Dedekind zeta values.
Example 3 expresses partial Eisenstein-Kronecker series associated
to an  imaginary quadratic ring as multiple Dedekind zeta values (see \cite{G2}, section 8.1).
Example 4 considers multiple Eisenstein-Kronecker series (for an alternative definition see \cite{G2}, 
Section 8.2). 
Examples 5 give the simplest multiple Dedekind zeta functions. 
Example 6 is a double Dedekind zeta function. 

In Section \ref{sec analytic}, we prove an analytic continuation of multiple Dedekind zeta functions, which allows us to consider special values of multiple Eisenstein series, examined by Gangl, Kaneko and Zagier (see \cite{ZGK}), as values of multiple Dedekind zeta functions, (see Examples 7, 8, 9 in Subsection \ref{subsec Eis}).
Examples 10 and 11 are particular cases of analytic continuation and of a multiple residue at 
$(1,\dots,1)$. The proof of analytic continuation is based on a generalization of Example 11 and a Theorem of Gelfand-Shilov (Theorem \ref{thm Gelfand-Shilov}).
At the end of Section \ref{sec remarks}, based on Examples 10 and 11, we state two conjectures about MDZV.

\section{Examples}
\label{sec examples}
We are going to present several examples of Riemann zeta values and multiple zeta values in order to introduce key examples of multiple Dedekind zeta value as iterated integrals. Instead of considering the ring of integers in a general number field, which we will do in the later sections, we will examine only the ring of Gaussian integers. Also, here we will ignore questions about convergence. Such questions will be addressed in Subsection \ref{sec cone}.

\subsection{Classical cases}
\label{subsec examples polylog}
Let us recall the $m$-th polylogarithm and its relation to Riemann zeta values.

If the first polylogarithm is defined as 
$$Li_1(x_1)=\int_0^{x_1}\frac{dx_0}{1-x_0}=\int_0^{x_1}(1+x_0+x_0^2+\dots)dx_0=
x_1+\frac{x_1^2}{2}+\frac{x^3_1}{3}+\dots$$
and  the second polylogarithm is
$$Li_2(x_2)=\int_0^{x_2}Li_1(x_1)\frac{dx_1}{x_1}
=
x_2+\frac{x_2^2}{2^2}+\frac{x^3_2}{3^2}+\dots$$
(Note that $\zeta(2)=Li_2(1)$),
then the $m$-th polylogarithm is defined  by iteration
\begin{equation}
\label{eq n-log}
Li_m(x_m)=\int_0^{x_m}Li_{m-1}(x_{m-1})\frac{dx_{m-1}}{x_{m-1}}.
\end{equation}
This is a presentation of the $m$-th polylogarithm as an iterated integral.
By a direct computation it follows that
$$Li_m(x)=x+\frac{x^2}{2^m}+\frac{x^3}{3^m}+\dots$$
and the relation 
$$\zeta(m)=Li_m(1)$$
is straightforward. 
Using Equation \ref{eq n-log}, we can express the $m$-th polylogarithm as
$$Li_m(x_m)=\int_{0<x_0<x_1<\dots<x_m}
\frac{dx_0}{1-x_0}\wedge\frac{dx_1}{x_1}\wedge\cdots\wedge\frac{dx_{m-1}}{x_{m-1}}.$$

Let $x_i=e^{-t_i}$. Then the $m$-th polylogarithm can be written in the variables $t_0,\dots,t_m$ in the following way
\begin{equation}
\label{eq polylog int}
Li_m(e^{-{t_m}})=\int_{t_0>t_1>\dots>t_m}\frac{dt_0\wedge\dots\wedge dt_{m-1}}{e^{t_0}-1}.
\end{equation}
This is achieved, first, by changing the variables in the differential forms
$$\frac{dx_0}{1-x_0}=\frac{d(-t_0)}{e^{t_0}-1},\text{ and }\frac{dx_i}{x_i}=d(-t_i),$$
and second, by reversing the bounds of integration $0<x_0<x_1<\dots<x_m$ v.s. 
$t_0>t_1>\dots>t_m$, which absorbs the sign.
As an infinite sum, we have
\begin{equation}
\label{eq polylog exp}
Li_m(e^{-t})=\sum_{n>0}\frac{e^{-nt}}{n^m}.
\end{equation}

In Subsection \ref{subsec examples polylog Gauss int}, we present a key analogy of Equations \eqref{eq polylog int} and \eqref{eq polylog exp} leading to Dedekind polylogarithms over the Gaussian integers.
Equations \eqref{eq polylog int} and \eqref{eq polylog exp} will be generalized to Dedekind polylogarithms in Subsection \ref{subsec Dpolylog} and to multiple Dedekind zeta values in Subsection \ref{subsec MDZV}. 

Below we present similar formulas for multiple 
polylogarithms with exponential variables. We will construct  their generalizations in Subsection \ref{subsec examples polylog Gauss int}.

Let us recall the definition of double logarithm
\begin{align*}Li_{1,1}(1,x_2)
&=
\int_0^{x_2}Li_1(x_1)\frac{dx_1}{1-x_1}
=
\int_0^{x_2}
\left(\sum_{n_1=1}^{\infty}\frac{x_1^{n_1}}{n_1}\right)
\left(\sum_{n_2=1}^{\infty}x_1^{n_2}\right)
\frac{dx_1}{x_1}
=\\
&=
\sum_{n_1,n_2=1}^{\infty}\frac{x_2^{n_1+n_2}}{n_1(n_1+n_2)}.
\end{align*}
Let $x_i=e^{-t_i}$. Then the $Li_{1,1}(1,e^{-t_2})$ can be written as an iterated integral in terms of the variables $t_0,t_1,t_2$ in the following way
$$Li_{1,1}(1,e^{-t_2})=\int_{t_0>t_1>t_2>0}\frac{dt_0\wedge dt_1}{(e^{t_0}-1)(e^{t_1}-1)}.$$
As an infinite sum, we have
\begin{equation}
\label{eq double polylog exp}
Li_{1,1}(1,e^{-t})=\sum_{n_1,n_2=1}^\infty\frac{e^{-(n_1+n_2)t}}{n_1(n_1+n_2)}.
\end{equation}

An example of a multiple zeta value is 
$$\zeta(1,2)=\sum_{n_1,n_2=1}^{\infty}\frac{1}{n_1(n_1+n_2)^2}=\int_0^1 Li_{1,1}(x_2)\frac{dx_2}{x_2}.$$

Thus, an integral representation of $\zeta(1,2)$ is
\begin{equation}
\label{eq z12}
\zeta(1,2)
=
\int_{t_0>t_1>t_2>0}\frac{dt_0}{(e^{t_0}-1)}\wedge \frac{dt_1}{(e^{t_1}-1)}\wedge dt_2.
\end{equation}

Similarly,
\begin{equation}
\label{eq z22}
\zeta(2,2)
=
\int_{t_0>t_1>t_2>t_3>0}\frac{dt_0}{(e^{t_0}-1)}\wedge dt_1\wedge\frac{dt_2}{(e^{t_2}-1)}\wedge dt_3.
\end{equation}

\subsection{Dedekind polylogarithms over the Gaussian integers}
\label{subsec examples polylog Gauss int}

In this Subsection, we are going to construct analogues of polylogarithms (and of some multiple polylogarithms), which we call Dedekind (multiple) polyologarithms over the Gaussian integers. 
We will denote by $f_{m}$ the $m$-th Dedekind polylogarithm, which will be an analogue the $m$-th polylogarithm $Li_m(e^{-t})$ with an exponential variable. Each of the analogues will have an integral representation, resembling an iterated integral and an infinite sum representation, resembling the classical Dedekind zeta values over the Gaussian integers. We also draw diagrams that represent  integrals in order to give a geometric view of the iterated integrals on membranes in dimension $2$. We will give examples of multiple Dedekind zeta values (MDZV) over the Gaussian integers, using the Dedekind (multiple) polylogarithms.

We are going to generalize Equations \eqref{eq polylog exp} and \eqref{eq double polylog exp} for (multiple) polylogarithms to their analogue over the Gaussian integers.  We will recall some properties and definitions related to Gaussian integers. For more information one may consider \cite{IR}.

By Gaussian integers we mean all numbers of the form $a+ib$, where $a$ and $b$ are integers and $i=\sqrt{-1}$. The ring of Gaussian integers is denoted by $\Z[i]$. We call the following set $C$ a {\it{cone}}
$$C=\N\{1+i,1-i\}=\{\alpha\in\Z[i]\mbox{ }|\mbox{ }\alpha=a(1+i)+b(1-i);\mbox{ } a,b \in \N\},$$
where $\N$ denotes the positive integers.
Note that $0$ does not belong to the cone $C$, since the coefficients $a$ and $b$ are positive integers.
We are going to use two sequences of inequalities
$$t_1>u_1>v_1>w_1\mbox{ and }t_2>u_2>v_2>w_2,$$
when we deal with a small number of iterations. The reason for introducing them is to make the examples easier to follow.
However, for generalizations to higher order of iteration we will use the following notation for the two 
sequences
$$t_{1,0}>t_{1,1}>t_{1,2}>t_{1,3}\mbox{ and }t_{2,0}>t_{2,1}>t_{2,2}>t_{2,3}.$$
We are going to define a function $f_1$, which will be an analogue of $Li_1(e^{-t})$. 
Let 
\begin{equation}
\label{eq Gauss f0}
f_0(C;t_1,t_2)=\sum_{\alpha\in C}\exp(-\alpha t_1 - \overline{\alpha}t_2).
\end{equation}
$$f_1(C,u_1,u_2)=\int^{u_1}_{\infty}\int^{u_2}_{\infty}f_0(C;t_1,t_2)dt_1\wedge dt_2.$$
We can draw the following diagram for the integral representing $f_1$.
\begin{center}
\begin{tikzpicture}
\draw[step=2cm] (0,0) grid (2,2);
\draw (0,-.5)node{$+\infty$};
\draw (1,-.5)node{$t_1$};
\draw (2,-.5)node{$u_1$};
\draw (-.5,0)node{$+\infty$};
\draw (-0.5,1)node{$t_2$};
\draw (-.5,2)node{$u_2$};
\draw (1,1)node{$f_0dt_1\wedge dt_2$};
\end{tikzpicture}
\end{center}
The diagram represents that the integrant is $f_0(C;t_1,t_2)dt_1\wedge dt_2$, depending on the variables $t_1$ and $t_2$, subject to the restrictions $+\infty> t_1>u_1$ and $+\infty>t_2>u_2$.

We need the following:
\begin{lemma}
(a) $$\int^u_{\infty} e^{-kt}dt=\frac{e^{-ku}}{k};$$

(b) Let $N(\alpha)=\alpha\overline{\alpha}.$ Then

$$\int^{u_1}_{\infty}\int^{u_2}_{\infty}\exp(-\alpha t_1 - \overline{\alpha}t_2)dt_1\wedge dt_2=
\frac{\exp(-\alpha u_1 - \overline{\alpha}u_2)}{N(\alpha)}.$$
\end{lemma}
The proof is straight forward.

Using the above Lemma, we obtain
$$f_1(C;u_1,u_2)
=
\sum_{\alpha\in C}\frac{\exp(-\alpha u_1- \overline{\alpha}u_2)}{N(\alpha)}
$$
We define a Dedekind dilogarithm $f_2$ by
\begin{align}
\nonumber
f_2(C;v_1,v_2)&=\int^{v_1}_{\infty}\int^{v_2}_{\infty}f_{1}(C;u_1,u_2)d u_1\wedge du_2=\\
\label{eq Ddilog}&=\int_{t_1>u_1>v_1;\mbox{ }t_1>u_1>v_2}f_0(C;t_1,t_2)d t_1\wedge dt_2\wedge d u_1\wedge du_2
\end{align}
We can associate a diagram to the integral representation of the Dedekind dilogarithm $f_2$ (see Equation \eqref{eq Ddilog}).

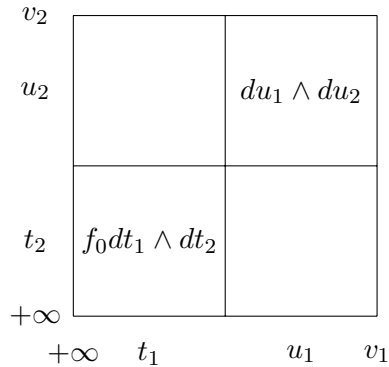
\begin{center}
\begin{tikzpicture}
\draw[step=2cm] (0,0) grid (4,4);
\draw (0,-.5)node{$+\infty$};
\draw (1,-.5)node{$t_1$};
\draw (3,-.5)node{$u_1$};
\draw (4,-.5)node{$v_1$};
\draw (-.5,0)node{$+\infty$};
\draw (-0.5,1)node{$t_2$};
\draw (-.5,3)node{$u_2$};
\draw (-.5,4)node{$v_2$};
\draw (1,1)node{$f_0dt_1\wedge dt_2$};
\draw (3,3)node{$du_1\wedge du_2$};
\end{tikzpicture}
\end{center}
The diagram represents that the variables under the integral are $t_1$, $t_2$, $u_1$, $u_2$, 
subject to the conditions $+\infty>t_1>u_1>v_1$ and $+\infty>t_2>u_2>v_2$. Also, the function $f_0$ in the diagram depends on the variables $t_1$ and $t_2$.

Similarly to Equation \eqref{eq n-log}, we define inductively the $m$-th Dedekind polylogarithm over the Gaussian integers
\begin{equation}
\label{eq X}
f_m(C;t_{1,m},t_{2,m})=\int^{t_{1,m}}_{\infty}\int^{t_{2,m}}_{\infty}f_{m-1}(C;t_{1,m-1},t_{2,m-1})dt_{1,m-1}\wedge dt_{2,m-1},
\end{equation}
where
$t_{1,0}>t_{2,1}>\dots>t_{1,m-1}>t_{1,m}\mbox{ and }t_{2,0}>t_{2,1}>\dots>t_{2,m-1}>t_{2,m}.$
The above integral is the key example of an iterated integral over a membrane, which is the topic of Subsection \ref{sec membrane}.

From Equation \eqref{eq X}, we can derive an analogue of the infinite sum representation of a polylogarithm
(see Equation  \eqref{eq polylog exp}).
\begin{equation}
\label{eq Y}
f_m(C;t_{1,m},t_{2,m})
=
\sum_{\alpha\in C}\frac{\exp(-\alpha t_{1,m} - \overline{\alpha}t_{2,m})}{N(\alpha)^m}.
\end{equation}
The above Equation gives an infinite sum representation of the $m$-th Dedekind polylogarithm over the Gaussian integers.

We derive one relation between the Dedekind $m$-polylogarithm $f_m$, a Dedekind zeta value over the Gaussian integers and a Riemann zeta value.
For arithmetic over the Gaussian integers one can consider \cite{IR}.
\begin{lemma} For the Dedekind polylogarithm $f_m$, associated to the above cone $C$, we have
\[f_m(C;0,0)=2^{-m}(\zeta_{\Q(i)}(m)-\zeta(2m)),\]
where $\zeta_{\Q(i)}(m)$ is a Dedekind zeta  value and $\zeta(2m)$ is a Riemann zeta value.
\end{lemma}
\proof
We are going to prove the following equalities, which give the lemma.
\begin{align}
f_m(C;0,0)&=
\sum_{\alpha\in C}\frac{1}{N(\alpha)^m}
=
2^{-m}\left(\sum_{(\alpha)\neq (0)\subset \Z[i]}-\sum_{\alpha\in \N}\right)\frac{1}{N((\alpha))^m}=\nonumber\\
&=2^{-m}(\zeta_{\Q(i)}(m)-\zeta(2m)),
\end{align}
The first equality follows from \eqref{eq Y}. The second and the third equalities relate our integral to classical zeta values.
The second equality uses two facts: (1) for the Gaussian integers the norm of an element $\alpha$,
 $N(\alpha)$, is equal to the norm of the principal ideal generated by $\alpha$, denoted by $N((\alpha))$, namely $N(\alpha)=N((\alpha))$.
 Recall that for the Gaussian integers the norm of an element $\alpha$,
 is $N(\alpha)=\alpha\overline{\alpha}$, and  the norm of a principal ideal $N((\alpha))$
  is equal to the number of elements in the quotient module
 \[N((\alpha))=\#|\Z[i]/(\alpha)|,\]  
  where
 \[(\alpha)=\alpha\Z[i]=\{\mu\in\Z[i]\mbox{ }|\mbox{ }\mu=\alpha\beta\text{ for some }\beta\in\Z[i]\}\] 
is view as a $\Z[i]$-submodule of $\Z[i]$.
(2) the set of non-zero principal ideals can be parametrized by the non-zero integers modulo the units. Since the units are $\pm 1,\pm i$, we have that $(\alpha)\subset \Z[i]$,  $(\alpha)\neq(0)$ can be parametrized by elements of the Gaussian integers with positive real part and non-negative imaginary part, which we will denote by $C_0$. Multiplying each element of $C_0$ by $1-i$, we obtain the union of the cone $C$ and the set $\{a+ai\mbox{ }|\mbox{ }a\in\N\}.$ Summing over $C_0$ gives the Dedekind zeta value. Summing over $(1-i)C$ gives $2{-m}\zeta_{\Q(i)}(m)$. Such a sum can be separated to a sum over $C$, which contributes $f_m$ and a sum over the set $\{a+ai\mbox{ }|\mbox{ }a\in\N\},$ which gives $2^{-m}\zeta(2m)$.\qed

Now we can define an analogue of the double logarithm $Li_{1,1}(1,e^{-t})$ over the  Gaussian integers, using the following integral representation
$$f_{1,1}(C;v_1,v_2)=\int^{v_1}_{\infty}\int^{v_2}_{\infty}
f_1(C;u_1,u_2)f_0(C;u_1,u_2)du_1\wedge du_2,$$
called a Dedekind double logarithm.
Such an integral will be considered as an example of an iterated integral over a membrane in Subsection 2.1. As an analog for Equation \eqref{eq Ddilog}, we can express $f_{1,1}$ only in terms of $f_0$ by
\[f_{1,1}(C;v_1,v_2)=\int_{t_1>u_1>v_1;\mbox{ }t_2>u_2>v_2}(f_0(C;t_1,t_2)d t_1\wedge dt_2)\wedge (f_0(u_1,u_2)d u_1\wedge du_2).\]
It allows us to associate a diagram to the Dedekind double logarithm $f_{1,1}$:
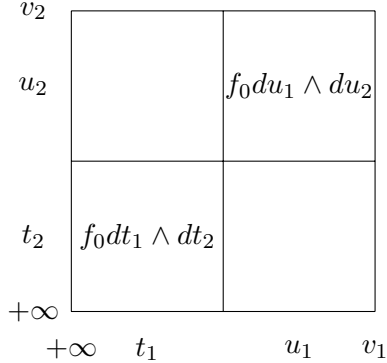
\begin{center}
\begin{tikzpicture}
\draw[step=2cm] (0,0) grid (4,4);
\draw (0,-.5)node{$+\infty$};
\draw (1,-.5)node{$t_1$};
\draw (3,-.5)node{$u_1$};
\draw (4,-.5)node{$v_1$};
\draw (-.5,0)node{$+\infty$};
\draw (-0.5,1)node{$t_2$};
\draw (-.5,3)node{$u_2$};
\draw (-.5,4)node{$v_2$};
\draw (1,1)node{$f_0dt_1\wedge dt_2$};
\draw (3,3)node{$f_0du_1\wedge du_2$};
\end{tikzpicture}
\end{center}
The variables $t_1$, $t_2$, $u_1$, $u_2$ in the diagram are variables in the integrant. They are subject to the conditions $t_1>u_1>v_1$ and $t_2>u_2>v_2$. Also, the lower left function $f_0$ in the diagram depends on the variables $t_1$ and $t_2$ and the upper right function $f_0$ depends on $u_1$ and $u_2$. 

The similarity between $f_{1,1}(C;v_1,v_2)$ and $Li_{1,1}(1,e^{-t_2})$ can be noticed by the infinite sum representation in the following:
\begin{lemma}
$$f_{1,1}(C;v_1,v_2)
=
\sum_{\alpha,\beta\in C}
\frac
{\exp(-(\alpha+\beta)v_1-(\overline{\alpha}+\overline{\beta})v_2)}
{N(\alpha)N(\alpha+\beta)}.$$
\end{lemma}
\proof
\begin{align*}
f_{1,1}(C;v_1,v_2)&=\int^{v_1}_{\infty}\int^{v_2}_{\infty}
f_1(C;u_1,u_2)f_0(C;u_1,u_2)du_1\wedge du_2=\\
&=\int^{v_1}_{\infty}\int^{v_2}_{\infty}
\sum_{\alpha\in C}
\frac
{\exp(-\alpha u_1-\overline{\alpha}u_2)}
{N(\alpha)}
\sum_{\beta\in C}
\exp(-\beta u_1-\overline{\beta}u_2)
du_1\wedge du_2
=\\
&=\int^{v_1}_{\infty}\int^{v_2}_{\infty}
\sum_{\alpha,\beta\in C}
\frac
{\exp(-(\alpha+\beta)u_1-(\overline{\alpha}+\overline{\beta})u_2)}
{N(\alpha)}
du_1\wedge du_2
=\\
&=\sum_{\alpha,\beta\in C}
\frac{\exp(-(\alpha+\beta)v_1-(\overline{\alpha}+\overline{\beta})v_2)}
{N(\alpha)N(\alpha+\beta)}.
\end{align*}
Similarly to the Dedekind double logarithm $f_{1,1}$, we define a multiple Dedekind polylogarithm
$$f_{1,2}(C,w_1,w_2)
=
\int^{w_1}_{\infty}\int^{w_2}_{\infty}
f_{1,1}(C;v_1,v_2)dv_1\wedge dv_2.$$
We can associate the following diagram to the multiple Dedekind polylogarithm $f_{1,2}$

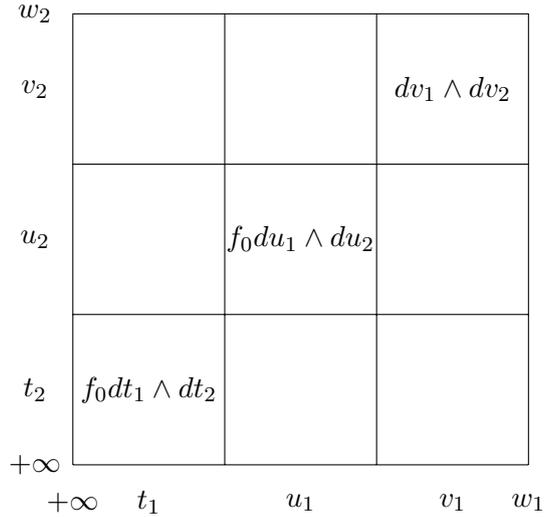
\begin{center}
\begin{tikzpicture}
\draw[step=2cm] (0,0) grid (6,6);
\draw (0,-.5)node{$+\infty$};
\draw (1,-.5)node{$t_1$};
\draw (3,-.5)node{$u_1$};
\draw (5,-.5)node{$v_1$};
\draw (6,-.5)node{$w_1$};
\draw (-0.5,0)node{$+\infty$};
\draw (-0.5,1)node{$t_2$};
\draw (-.5,3)node{$u_2$};
\draw (-.5,5)node{$v_2$};
\draw (-.5,6)node{$w_2$};
\draw (1,1)node{$f_0dt_1\wedge dt_2$};
\draw (3,3)node{$f_0du_1\wedge du_2$};
\draw (5,5)node{$dv_1\wedge dv_2$};
\end{tikzpicture}
\end{center}
The diagram represents the following: The variables of the integrant are $t_1$,  $t_2$, $u_1$, $u_2$, $v_1$,  $v_2$. 
The variables are subject to the conditions $t_1>u_1>v_1>w_1$ and $t_2>u_2>v_2>w_2$.
The lower left function $f_0$ depends on the variables $t_1$ and $t_2$. And the middle function $f_0$ depends on $u_1$ and $u_2$.  Thus, the diagram represents the following integral:
\begin{align}
\label{eq pre f12}
&f_{1,2}(C;w_1,w_2)=\\
\nonumber
&=\int_{D_{w_1,w_2}}
(f_0(C;t_1,t_2)dt_1\wedge dt_2)\wedge
(f_0(C;u_1,u_2)du_1\wedge du_2)\wedge (dv_1\wedge dv_2),
\end{align}
where the domain of integration is
\[D_{w_1,w_2}=\{(t_1,t_2,u_1,u_2,v_1,v_2)\in \R^6\mbox{ }|\mbox{ }t_1>u_1>v_1>w_1\mbox{ and } t_2>u_2>v_2>w_2\}\]
A direct computation leads to
$$f_{1,2}(C;w_1,w_2)
=\sum_{\alpha,\beta\in C}
\frac{\exp(-(\alpha+\beta)w_1-(\overline{\alpha}+\overline{\beta})w_2)}
{N(\alpha)N(\alpha+\beta)^2}.
$$
We define a multiple Dedekind zeta value as 
$$\zeta_{\Q(i);C}(1,1;2,2)
=
f_{1,2}(C;0,0)
=
\sum_{\alpha,\beta\in C}
\frac{1}{N(\alpha)N(\alpha+\beta)^2}.
$$

Now let us give a relation between multiple Dedekind zeta values and iterated integrals. 
We use the following pair of inequalities in the following Sections
$t_{1,0}>t_{1,1}>t_{1,2}\mbox{ and }t_{2,0}>t_{2,1}>t_{2,2},$
instead of
$t_1>u_1>v_1>w_1\mbox{ and }t_2>u_2>v_2>w_2,$
since using such notation it is easier to write higher order iterated integrals.
In this notation, from Equation \eqref{eq X} , we obtain
\begin{align}
\label{eq f2}&f_{2}(C;t_{1,2},t_{2,2})
=\\
\nonumber
&=
\int_{t_{1,0}>t_{1,1}>t_{1,2};\mbox{ } t_{2,0}>t_{2,1}>t_{2,2}}
(f_0(C;t_{1,0},t_{2,0})dt_{1,0}\wedge dt_{2,0})
\wedge 
(dt_{1,1}\wedge dt_{2,1}).
\end{align}
and
\begin{align}
\label{eq f11}
&f_{1,1}(C;t_{1,2},t_{2,2})
=\\
\nonumber
&=
\int_{t_{1,0}>t_{1,1}>t_{1,2};\mbox{ } t_{2,0}>t_{2,1}>t_{2,2}}
(f_0(C;t_{1,0},t_{2,0})dt_{1,0}\wedge dt_{2,0})
\wedge 
(f_0(C;t_{1,1},t_{2,1})dt_{1,1}\wedge dt_{2,1}).
\end{align}

In the next Section, we generalize the (iterated) integrals appearing in Equations \eqref{eq X}
\eqref{eq pre f12}, \eqref{eq f2}, and \eqref{eq f11}, called iterated integrals over a membrane, (see Definition \ref{def it int n-forms}). 

The next two Examples are needed in order to relate multiple Dedekind zeta values to values of Eisenstein series and values of multiple Eisenstein series (see \cite{ZGK},  see also Examples 7, 8, 9 at the end of Section \ref{sec MDZF}).
The integrals below will present another type of iterated integral on a membrane, (see Definition \ref{def it int n-forms and 1-forms}) leading to a multiple Dedekind zeta values. 

We define the following iterated integral to be a multiple Dedekind zeta value 
\begin{equation}
\label{eq z32}
\zeta_{\Q(i);C}(3;2)
=
\int_{t_1>u_1>v_1>0;\mbox{ }t_2>u_2>0}
f_0(C;t_1,t_2)
dt_1\wedge du_1\wedge dv_1\wedge dt_2\wedge du_2.
\end{equation}
The reason for such a definition is its infinite sum representation
\begin{equation}
\label{eq sum z32}
\zeta_{\Q(i);C}(3;2)=\sum_{\alpha\in C}\frac{1}{\alpha^3\overline{\alpha}^2},
\end{equation}
which can be achieved essentially in the same way as for the other multiple Dedekind zeta values.
We can associate the following diagram to the integral representation of $\zeta_{\Q(i);C}(3;2)$ in Equation \eqref{eq z32}.

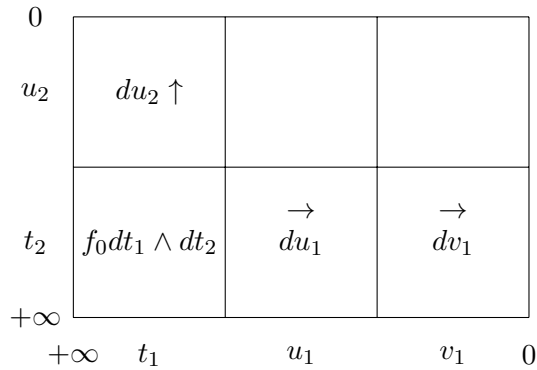
\begin{center}
\begin{tikzpicture}
\draw[step=2cm] (0,0) grid (6,4);
\draw (0,-.5)node{$+\infty$};
\draw (1,-.5)node{$t_1$};
\draw (3,-.5)node{$u_1$};
\draw (5,-.5)node{$v_1$};
\draw (6,-.5)node{$0$};
\draw (-0.5,0)node{$+\infty$};
\draw (-0.5,1)node{$t_2$};
\draw (-.5,3)node{$u_2$};
\draw (-.5,4)node{$0$};
\draw (1,1)node{$f_0dt_1\wedge dt_2$};
\draw (1,3)node{$du_2\uparrow$};
\draw (3,1.4)node{$\rightarrow$};
\draw (3,1)node{$du_1$};
\draw (5,1.4)node{$\rightarrow$};
\draw (5,1)node{$dv_1$};
\end{tikzpicture}
\end{center}
The arrows in the diagram signify the direction of decrease of the variables in differential $1$-forms.
It is important to consider the variables only in the horizontal direction and then only in vertical direction. 
In horizontal direction, we have 
$f_0(C;t_1,t_2)dt_1$ followed by $du_1$ and $dv_1$, The integration with respect to the variables $t_1,u_1,v_1$ leads to $1/\alpha^3$ in the summation of Equation \eqref{eq sum z32}. In vertical direction, we have a double iteration. First we have $f_0(C;t_1,t_2)dt_2$ followed by $u_2$. That leads to $1/\overline{\alpha}^2$ in the summation in Equation \eqref{eq sum z32}.

Diagrams associated to the integral representation of $\zeta_{\Q(i);C}(3;2)$ (Equations \eqref{eq z32}) are not unique. Alternatively, we could have used the diagrams

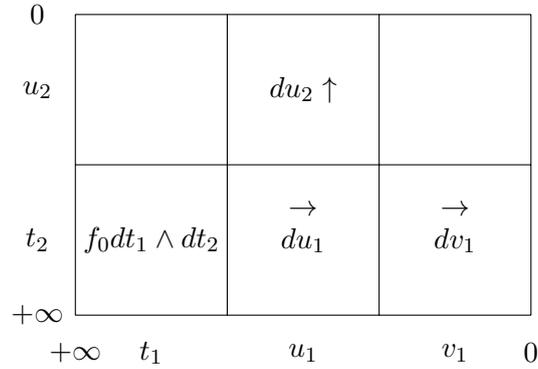
\begin{center}
\begin{tikzpicture}
\draw[step=2cm] (0,0) grid (6,4);
\draw (0,-.5)node{$+\infty$};
\draw (1,-.5)node{$t_1$};
\draw (3,-.5)node{$u_1$};
\draw (5,-.5)node{$v_1$};
\draw (6,-.5)node{$0$};
\draw (-0.5,0)node{$+\infty$};
\draw (-0.5,1)node{$t_2$};
\draw (-.5,3)node{$u_2$};
\draw (-.5,4)node{$0$};
\draw (1,1)node{$f_0dt_1\wedge dt_2$};
\draw (3,3)node{$du_2\uparrow$};
\draw (3,1.4)node{$\rightarrow$};
\draw (3,1)node{${du_1}$};
\draw (5,1.4)node{$\rightarrow$};
\draw (5,1)node{$dv_1$};
\end{tikzpicture}
\end{center}
or
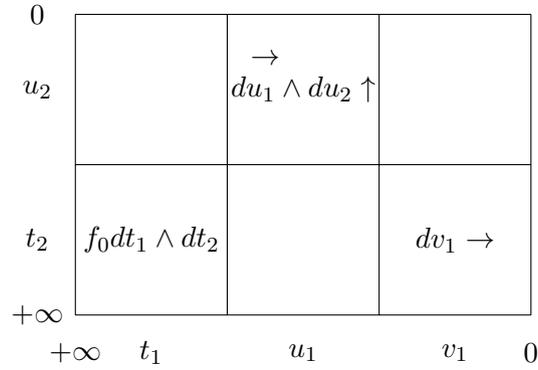
\begin{center}
\begin{tikzpicture}
\draw[step=2cm] (0,0) grid (6,4);
\draw (0,-.5)node{$+\infty$};
\draw (1,-.5)node{$t_1$};
\draw (3,-.5)node{$u_1$};
\draw (5,-.5)node{$v_1$};
\draw (6,-.5)node{$0$};
\draw (-0.5,0)node{$+\infty$};
\draw (-0.5,1)node{$t_2$};
\draw (-.5,3)node{$u_2$};
\draw (-.5,4)node{$0$};
\draw (1,1)node{$f_0dt_1\wedge dt_2$};
\draw (2.5,3.4)node{$\rightarrow$};
\draw (3,3)node{$du_1\wedge du_2\uparrow$};
\draw (5,1)node{$dv_1\rightarrow$};
\end{tikzpicture}
\end{center}

Consider the following diagram, associated to a more complicated MDZV, for the purpose of establishing notation.

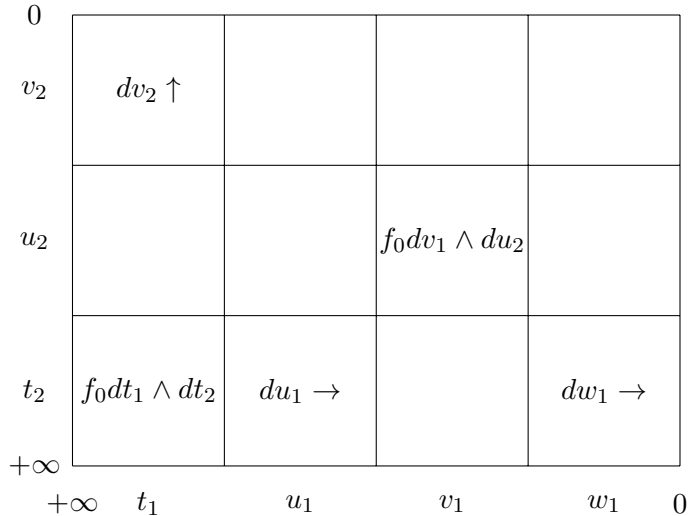
\begin{center}
\begin{tikzpicture}
\draw[step=2cm] (0,0) grid (8,6);
\draw (0,-.5)node{$+\infty$};
\draw (1,-.5)node{$t_1$};
\draw (3,-.5)node{$u_1$};
\draw (5,-.5)node{$v_1$};
\draw (7,-.5)node{$w_1$};
\draw (8,-.5)node{$0$};
\draw (-0.5,0)node{$+\infty$};
\draw (-0.5,1)node{$t_2$};
\draw (-.5,3)node{$u_2$};
\draw (-.5,5)node{$v_2$};
\draw (-.5,6)node{$0$};
\draw (1,1)node{$f_0dt_1\wedge dt_2$};
\draw (1,5)node{$dv_2\uparrow$};
\draw (3,1)node{$du_1\rightarrow$};
\draw (5,3)node{$f_0dv_1\wedge du_2$};
\draw (7,1)node{$dw_1\rightarrow$};
\end{tikzpicture}
\end{center}
The diagram encodes that $+\infty>t_1>u_1>v_1>w_1>0$ and $+\infty>t_2>u_2>v_2>0$.
Consider the horizontal direction of the diagram and Equation \eqref{eq z22}. We have $f_0dt_1$, followed by $du_1$, $f_0dv_1$ and $dw_1$. That gives an analogue of $\zeta(a,c)=\zeta(2,2)$ in horizontal direction for $(a,c)=(2,2)$.
Consider the vertical direction of the diagram and Equation \eqref{eq z12}.  We have $f_0dt_2$, followed by  $f_0du_2$ and $dv_2$. That  gives an analogue of $\zeta(b,d)=\zeta(1,2)$ in vertical direction, for $(b,d)=(1,2)$. We write
\[\zeta_{\Q(i);C,C}(a,b;c,d)=\zeta_{\Q(i);C,C}(2,1;2,2)\]
for the multiple Dedekind zeta function associated to the above diagram.
We leave proof of the following statement is left to the reader
 \[\zeta_{\Q(i);C,C}(2,1;2,2)=\sum_{\alpha,\beta\in C}\frac{1}{\alpha^2\overline{\alpha}^1
 (\alpha+\beta)^2(\overline{\alpha}+\overline{\beta})^2}.\]
 
In Section \ref{sec MDZF}, 
we use iterated integrals over a membrane 
to define multiple Dedekind zeta values associated to any number field. 

\section{Arithmetic and geometric tools}
\subsection{Iterated integrals on a membrane}
\label{sec membrane}  
Let $D$ be a domain defined  in terms of the real variables $t_{i,j}$ for $i=1,\dots,n$ and $j=1,\dots,m$, by
$$D=\{(t_{1,1},\dots,t_{n,m})\in \R^{nm}\mbox{ }|\mbox{ } t_{i,1}>t_{i,2}>\dots>t_{i,m}>0\text{ for }i=1,\dots,n\}.$$
For each $j=1,\dots,m$, let $\omega_j$ be a differential $n$-form on $\C^n$.
Let
$$g:(0,+\infty)^n\rightarrow \C^n$$
be a smooth map, whose pull-back sends the coordinate-wise foliation on $\C^n$ 
to a coordinate-wise foliation on $(0,+\infty)^n$.
We will  call such a map {\it{a membrane}}.
One should think of the $n$-forms $g^*\omega_j$ as an 
analogue of $f_0(C;t_1,t_2)dt_1\wedge dt_2$
from Equation \eqref{eq Gauss f0}.
\begin{definition}
\label{def it int n-forms} An iterated integral on a membrane $g$,
in terms of $n$-forms $\omega_j$, $j=1,\dots,m$,  is defined as
\begin{equation}
\label{eq it int n-forms}
\int_{g}\omega_1\dots\omega_m
=
\int_D\bigwedge_{j=1}^m g^*\omega_j(t_{1,j},\dots,t_{n,j}).
\end{equation}
\end{definition}

\begin{definition}
\label{def shuffle of sets}
A shuffle between two ordered sets
$$S_1=\{1,\dots,p\}$$
and
$$S_2=\{p+1,\dots,p+q\}$$
is a permutation $\tau$ of the union $S_1\cup S_2$, such that
\begin{enumerate}

\item
for $a,b\in S_1$, we have $\tau(a)<\tau(b)$ if $a<b$;
\item
for $a,b\in S_2$, we have $\tau(a)<\tau(b)$ if $a<b$;
\end{enumerate}
We denote the set of shuffles between two ordered sets of orders
$p$ and $q$, respectively, by $Sh(p,q)$.
\end{definition}
The definition of an iterated integral on a membrane is associated
with the following objects:

\begin{enumerate}

\item 
$g:(0,+\infty)^n\rightarrow \C^n$,
a membrane (that is a smooth map, whose pull-back sends the coordinate-wise foliation on $\C^n$ to a coordinate-wise foliation on $(0,+\infty)^n$).


\item $\omega_1,\ldots, \omega_m$ differential $n$-forms on $\C^n$;

\item $m_i$ copies of differential $1$-forms $\di z_i$ on $\C^n$, for $i=1,\dots,n$;

\item  a shuffle $\tau_i \in Sh(m,m_i)$ for each $i=1,\dots,n$;

\item $\tau=(\tau_1,\dots,\tau_n)$, the set of $n$ shuffles $\tau_1,\dots,\tau_n.$

\end{enumerate}

\begin{definition}
\label{def it int n-forms and 1-forms} Given the above data, we
define an iterated integral on a membrane $g$, involving $n$-forms
and $1$-forms, as
\begin{align}
&\nonumber \int_{(g,\tau)}\omega_1\dots\omega_m(\di z_1)^{m_1}\dots(\di z_n)^{m_n}
=\\
\label{eq def MDZV}=
&\int_D
\left(\prod_{j=1}^mg^*\omega_j(t_{1,\tau_1(j)},\dots,t_{n,\tau_n(j)})\right)
\bigwedge_{i=1}^n\bigwedge_{j=1}^{m+m_i} g^*dz_{i,j},
\end{align}
where $t_{i,j}=g^*z_{i,j}$ and also $t_{i,j}$ belong to the domain
$$D=\{(t_{1,1},\dots,t_{n,m})\in \R^{mn}\mbox{ }|\mbox{ }t_{i,1}>t_{i,2}>\dots>t_{i,m+m_i}>0\}.$$
\end{definition}

{\bf{Remark:}} Comparing the Definitions \ref{def it int n-forms} and \ref{def it int n-forms and 1-forms}, one can notice that there is no sign occurring. The reason for that is the following: 
\begin{enumerate}
\item In Definition \ref{def it int n-forms} we use a domain $D$, whose coordinates are ordered  by $t_{1,1},\dots,t_{n,1}$, $t_{1,2},\dots,t_{n,2}$, $\dots,t_{1,m},\dots,t_{n,m}.$ It is the same as the order of the differential 1-forms under the integral in Equation \eqref{eq it int n-forms}. 
\item  In Definition \ref{def it int n-forms and 1-forms}
 we use a domain $D$, whose coordinates are ordered  by $t_{1,1},\dots,t_{1,m+m_1},$ $t_{2,1},\dots,t_{2,m+m_2}$, $\dots,t_{n,1},\dots,t_{n,m+m_n}.$ It is the same as the order of the differential 1-forms under the integral in Equation \eqref{eq def MDZV}.
\end{enumerate}
Thus, if $m_1=\cdots =m_n=0$, both definitions lead to the same value, since the permutation of the differential forms coincides with the permutation of the coordinates of the domain of integration. Thus, the change of orientation of the domain of integration coincides with the sign of  permutation acting on the differential forms.

\begin{theorem}(homotopy invariance)
\label{thm homotopy invariance}
The iterated integrals on membranes from Definition \ref{def it int n-forms and 1-forms}
are homotopy invariant, when the homotopy preserves the boundary of the membrane.
\end{theorem}
\proof
Let $g$ be a homotopy between the two membranes $g_0$ and $g_1$. Let 
\[\Omega=
\left(\prod_{j=1}^m\omega_j(z_{1,\tau_1(j)},\dots,z_{n,\tau_n(j)})\right)
\bigwedge_{i=1}^n\bigwedge_{j=1}^{m+m_i} dz_{i,j}\]

Note that $\Omega$ is a closed form, 
since $\omega_i$ is a form of top dimension and since $dz_{i,j}$ is closed.
By Stokes Theorem, we have
\begin{align}
\nonumber 0&=\int_{s=0}^{s=1}\int_D g^*d\Omega=\\
\label{eq g1 g0}
&=\int_{(g_1,\tau)}\Omega-\int_{(g_0,\tau)}\Omega\pm\\
\label{eq boundary1}
&\pm\int_{s=0}^{s=1}\sum_{i=1}^n\sum_{j=1}^{m+m_i}\int_{D|(z_{i,j}=z_{i,j+1})}g^*\Omega\pm\\
\label{eq boundary2}
&\pm\int_{s=0}^{s=1}\sum_{i=1}^n\int_{D|(z_{i,m+m_i}=0)}g^*\Omega
\end{align} 
We want to show that the difference in the terms in \eqref{eq g1 g0} is zero. It is enough to show that 
each of the terms \eqref{eq boundary1} and \eqref{eq boundary2} are zero. 
If $z_{i,j}=z_{i,j+1}$, then the wedge of the corresponding differential forms will vanish. Thus the terms in \eqref{eq boundary1} are zero. If $z_{i,m+m_i}=0$ then $dt_{i,m+m_i}=0$, defined via the pull-back $g^*$. Then  the terms \eqref{eq boundary2} are equal to zero.\qed

\subsection{Cones and geometric series}
\label{sec cone}
Let $n=[K:\Q]$ be the degree of the
number field $K$ over $\Q$. Let ${\cal{O}}_K$ be the ring of
integers in $K$. And let $U_k$ be the group of units in $K$. We
are going to use an idea of Shintani \cite{C} by examining
Dedekind zeta functions in terms of a cone inside the ring of
integers.

We define cone for any number ring. The meaning of  cones is roughly the following:  summation over the elements of  finitely many cones
would give multiple Dedekind zeta values or multiple Dedekind zeta functions.

\begin{definition}
We define a {\it{cone}} $C$ to be
\[C=\N\{e_1,\dots,e_k\}=\{\alpha \in {\cal{O}}_K \mbox{ }|\mbox{ }\alpha=a_1e_1+\cdots a_ke_k\mbox{ for }e_i\in {\cal{O}}_K\mbox{ and }a_i\in \N\}\]
with generators $e_1,\dots,e_k$.
\end{definition}

For the next definition, we are going to use that a number field $K$ can be viewed as an $n$-dimensional vector space over the rational numbers $\Q$.

\begin{definition}
An {\it{unimodular cone}} is a cone with generators $e_1,\dots,e_k$ such that $e_1,\dots,e_k$ as elements of $K$ are linearly independent over $\Q$, when we view the field $K$ as a vector space over $\Q$.
\end{definition}

Note that if $C$ is an unimodular cone then $0\notin C$, since $e_1,\dots,e_k$ are linearly independent over $\Q$ and the coefficients $a_1,\dots,a_n$ are positive integers.

\begin{definition}
We call $C$ an {\it{unimodular simple cone}} if for any embedding $\sigma_i$  of $K$ into the complex numbers, $\sigma_i:K\rightarrow \C$  and a suitable branch of the functions $\arg(z)$,
we have that the closure of the set $\arg(\sigma(C))$ is an interval $[\theta_0,\theta_1]$, such that its lengths is less than $\pi$, namely, $\theta_1-\theta_0\in[0,\pi)$.
\end{definition}

In particular, the cone 
$$C=\{\alpha\in \Z[i]\mbox{ }|\mbox{ }\alpha=a(1+i)+b(1+i), a,b\in \N\},$$ 
considered in Subsection \ref{subsec examples polylog Gauss int},
is an unimodular simple cone, since
$\arg(\sigma_1(\alpha))\in(-\pi/4,\pi/4)$
and
$\arg(\sigma_2(\alpha))\in(-\pi/4,\pi/4)$,
for each $\alpha\in C$. The maps $\sigma_1$ and $\sigma_2$ are complex conjugates of each other.

\begin{definition} (Dual cone)
\label{def dual cone}
For an unimodular simple cone $C$ with generators $e_1,\dots,e_k$, we define a dual cone of $C$ to be
\[C^*=\{(z_1,\dots,z_n)\in \C^n\mbox{ }|\mbox{ }Re(z_i\sigma_i(e_j))>0\mbox{ for }i=1,\dots,n,\mbox{ and }j=1\dots,k\}\]
\end{definition}
Clearly, if $C$ is an unimodular simple cone then the dual cone $C^*$ is a non-empty set. One can prove that by considering each coordinate of $C^*$, separately.  

\begin{definition} 
\label{def f0}
For an unimodular simple cone $C$, we define a function
\begin{equation}
f_0(C;z_1,\dots,z_n)=\sum_{\alpha\in C}\exp(-\sum_{i=1}^n\sigma_i(\alpha)z_i),
\end{equation}
where $\sigma_1,\dots,\sigma_n$ are all embeddings of the number field $K$ into the complex numbers $\C$ 
and the domain of the function $f_0$ is the dual cone $C^*$.
\end{definition}

\begin{lemma}
\label{lemma f0 convergence}
The function $f_0$ is uniformly convergent for $(z_1,\dots,z_n)$ in any compact subset $B$ of the dual cone $C^*$ of an unimodular simple cone $C$.
\end{lemma}

\proof From the Definition \ref{def dual cone}, we have $\Re(\sigma_i(e_j)z_i)>0$. Let
\begin{equation}
\label{eq yj}
y_j=\prod_{i=1}^n\exp(-\sigma_i(e_j)z_i).
\end{equation}
Then $|y_j|<1$ on the domain $B$. 
Moreover, $|y_i|$ achieves a maximum 
on the compact subset $B$.
Let $|y_j|\leq c_j<1$ on the compact set $B$ for some constant $c_j$.
Then the rate of convergence of the geometric sequence in $y_j$ 
is uniformly bounded by $c_j$ on the compact set $B$. 
Therefore, we have a uniform convergence of the geometric series in $y_j$.
The function $f_0(C;z_1,\dots,z_n)$ 
is a product of $k$ geometric series in the variables $y_1,\dots,y_k$ 
each of which is uniformly bounded in absolute value by the  constants  $c_1,\dots,c_k$ on the domain $B$, respectively. Then, we obtain that 
\begin{equation}
\label{eq f0}
f_0(C;z_1,\dots,z_n)=\prod_{j=1}^k\frac{y_j}{1-y_j}.
\end{equation}
\qed
\begin{corollary}
\label{cor cont f0} The function $f_0(C;z_1,\dots,z_n)$ has
analytic continuation to all values of $z_1,\dots,z_n$, except at
$$\sum_{i=1}^n\sigma_i(e_j)z_i\in 2 \pi i\Z,$$
for $j=1,\dots,k.$
\end{corollary}
\proof Using the product formula \ref{eq f0} in terms of geometric
series in $y_j$, we see that the right hand side of \ref{eq f0} makes sense
for all $y_j\neq 1$. This gives analytic continuation from the domain $C^*$ to the domain consisting of points $(y_1,\dots,y_n)$ with $y_i\neq 1$. 
\qed
\begin{definition} (positive cone)
We call $C$ a {\it{positive unimodular simple cone}} if $C$ is an unimodular simple cone 
and the product of the positive real coordinates is in $C^*$, namely
 \[(\R_{>0})^n\subset C^*\]
 as subsets of $\C^n$.
\end{definition}

\begin{lemma}
\label{lemma positive cone}
If $C$ is an unimodular simple cone then for some $\alpha\in {\cal{O}}_K$ we have that 
\[\alpha C=\{\alpha\beta\mbox{ }|\mbox{ }\beta\in C \}\] 
is a positive unimodular simple cone.
\end{lemma}
\proof In order to find such an elements $\alpha$, we need to recall properties of real or complex 
embeddings of a number field $K$. 

The degree of a number field $n=[K:\Q]$ is the dimension of $K$ as vector space over $\Q$. 
Then there are exactly $n$ distinct embeddings $K\rightarrow \C$. Let the first $r_1$ embeddings, 
$\sigma_1,\dots,\sigma_{r_1},$ be the ones whose image is inside the real numbers. They are called 
real embeddings. Let the next $r_2$ embeddings be complex embeddings, which are not pair-wise complex conjugates of each other. Let us denote them by $\sigma_{r_1+1},\dots,\sigma_{r_1+r_2}.$ Let the 
last $r_2$ embeddings be the complex conjugates of previously counted complex embeddings, 
namely,
\[\sigma_{r_1+r_2+i}(\beta)=\overline{\sigma_{r_1+i}(\beta)},\]
for $i=1,\dots,r_2$. We also have that $n=r_1+2r_2$.

Let $V_\R$ be a $n$ dimensional real vector subspace of $\C^n$ defined in the following way:
\begin{align*}
V_\R=\{(z_1,\dots,z_n)\in \C^n\mbox{ }|\mbox{ }&(z_1,\dots,z_{r_1})\in \R^{r_1},\mbox{  }\\ 
&(z_{r_1+1},\dots,z_{r_1+r_2})\in \C^{r_2}, \mbox{ and } \\
&(z_{r_1+r_2+1},\dots,z_{r_1+2r_2})=(\overline{z}_{r_1+1},\dots,\overline{z}_{r_1+r_2})\}
\end{align*}
Now, we proceed with the proof of the Lemma in six Steps.

Step 1. $K$ is dense in $V_\R$.

Step 2. $V_\R\cap C^*$ in non-empty. 

Step 3. $V_\R\cap C^*$ is an open subset of $V_\R$.

Step 4. $K\cap C^*$ is non-empty.

Step 5. ${\cal{O}}_K\cap C^*$ is non-empty.

Step 6. $\alpha C$ is a positive unimodular simple cone for any $\alpha\in {\cal{O}}_K\cap C^*$.

 (Step 1)
Recall also the product of the $n$ embeddings of $K$ to the space $V_\R$,
\[\prod_{i=1}^n\sigma_i:K\rightarrow V_\R,\]
mapping $\beta\in K$ to $(\sigma_1(\beta),\dots,\sigma_n(\beta))\in V_\R$
has a dense image. 

(Step 2)
Indeed, let $z_i$ be the $i$-th coordinate of $C^*$.
The first $r_1$ coordinates $z_1,\dots,z_{r_1}$ of $C^*$ can be real numbers (positive or negative), since $\sigma_i(K)\subset\R$ for $i=1,\dots,r_1$. Thus, the first $r_i$ coordinates can be both in $V_\R$ and in $C^*$. For the coordinates $z_{r_1+1},\dots,z_{r_1+r_2}$ of $C^*$ there are no restrictions when we intersect $C^*$ with $V_\R$. For the last $r_2$ coordinates of $C^*$ we must  have that $z_{r_1+r_2+i}=\overline{z}_{r_1+i}$ in order 
for the coordinates to be in the intersection $C^*\cap V_\R$.
Since, 
$\sigma_{r_1+r_2+i}(\beta)=\overline{\sigma_{r_1+i}(\beta)}$, we have the conditions on
the $(r_1+i)$-coordinate and on the $(r_1+r_2+i)$-coordinate of a point in $C^*$ to be in $V_\R$ are
\[Re(z_{r_1+i}\sigma_{r_1+i}(\beta))>0,\]
for $\beta\in C$ and $z_{r_1+r_2+i}=\overline{z}_{r_1+i}$. The last condition implies that
\[Re(z_{r_1+r_2+i}\sigma_{r_1+r_2i}(\beta))=Re(\overline{z_{r_1+i}\sigma_{r_1+i}(\beta)})>0.\]
Thus, such a point $(z_1,\dots,z_n)$ is in $C^*\cap V_\R$.

 (Step 3)
It is true, since $C^*$ is an open subsets of $\C^n$

 (Step 4) Since $K$ is dense in $V_\R$ (Step 1) and $V_\R\cap C^*$ is open in $V_\R$ (Steps 2 and 3), we have that $K\cap C^*$ is non-empty.

 (Step 5) If $\alpha\in K\cap C^*$ then for some positive integer $L$, we have that  $L\alpha\in {\cal{O}}_K$, and also, $L\alpha\in C^*$, since $C^*$ is invariant under rescaling by a positive (real) number $L$.

(Step 6) 
Let $(t_1,\dots,t_n)\in \R^n_{>0}$ and let $\beta\in C$. Put $z_i=t_i\sigma_i(\alpha)$. Then 
$(z_1,\dots,z_n)\in C^*$.
\[Re(t_i\sigma_i(\alpha\beta))=Re(t_i\sigma_i(\alpha)\sigma_i(\beta))=Re(z_i\sigma_i(\beta))>0.\]
\qed


\subsection{Cones and ideals}
\label{subsec DZ and cones}
In this Subsection, we are going to examine union of cones that give a fundamental domain of the ring of integers ${\cal{O}}_K$ modulo the group of units $U_K$. We also examine a fundamental domain of an ideal $\mathfrak{a}$ modulo the group of units $U_K$.  
\begin{definition}
\label{def fund domain}
We define $M$ as
a fundamental domain of
$${\cal{O}}_K-\{0\}\mbox{ }\mathrm{mod}\mbox{ }U_k.$$
For an ideal $\mathfrak{a}$, let
$$M(\mathfrak{a})=M\cap \mathfrak{a}.$$
\end{definition}
\begin{lemma}
\label{lemma fund domain}
\label{lemma fundamental domain} For any ideal $\mathfrak{a}$ the set 
$M(\mathfrak{a})$ can be written as a finite disjoint union of
unimodular simple cones.
\end{lemma}
\proof It is a simple observation that $M(\mathfrak{a})$ can
be written as a finite union of unimodular cones. We have to show
that we can subdivide each of the unimodular cones into finite
union of unimodular simple cones.

Let $\sigma_1,\dots,\sigma_{r_1}$ be the real embeddings of the number field $K$ and let $\sigma_{r_1+1},\dots,\sigma_{r_1+r_2}$ be the non-conjugate complex embeddings of $K$. 
We define 
\[T=\{-1,1\}^{r_1}\times (S^1)^{r_2}.\]
Let $C$ be an unimodular cone. We define a map $J$ by 
\begin{align*}
J:C&\rightarrow T\\
\alpha&\mapsto \left(\frac{\sigma_1(\alpha))}{|\sigma_1(\alpha))|},\cdots,\frac{\sigma_{r_1+r_2}(\alpha))}{|\sigma_{r_1+r_2}(\alpha))|}\right)
\end{align*}
Denote it by $\bar{C}$ the closure of the image of $J$ in $T$. Then one can cut the
cone $C$ into finitely many  cones  $C_i$ such that for $C_i$ and
any embedding $\sigma$ of $K$ into $\C$, we have that
$\arg(\sigma(C))\in[\theta_0,\theta_1]$, for $\theta_1-\theta_0\in [0,\pi)$. 
Then $C_{i}$ is an unimodular simple cone.
Thus, the cones
$C_{i}$'s are finitely many unimodular simple cones, whose (disjoint) union
gives the set $M(\mathfrak{a}).$
 \qed

For an element $\alpha$ in a ring of integers ${\cal{O}}_K$, 
denote by $(\alpha)$ the principal ideal generated by $\alpha$.
Then $N_{K/\Q}((\alpha))$ denoted the norm of the principal ideal generated by $\alpha$. 
We have that $N_{K/\Q}((\alpha))$ is a positive integer equal to the number of elements in the quotient  
${\cal{O}}_K/(\alpha)$. Also $N_{K/\Q}(\alpha)$ is the norm of the algebraic number $\alpha$. This is equal to the product of all of its Galois conjugates, which is an integer, possibly a negative integer. 
We always have that $N_{K/\Q}((\alpha))=|N_{K/\Q}(\alpha)|$.

However, for elements of an unimodular simple cone, we can say more.
\begin{lemma}
Let $C$ be an unimodular simple cone.
Then for every $\alpha\in C$, we have
$$N_{K/\Q}((\alpha))=\e(C)N_{K/\Q}(\alpha),$$
where $\e(C)=\pm 1$ depends only on the cone $C$, not on $\alpha$.
\end{lemma}
\proof Note that on the left we have a norm of an ideal and on the
right we have a norm of a number. Since $C$ is a simple cone, we
have that for all real embeddings $\sigma:K\rightarrow \R$, the
signs of $\sigma(\alpha)$ and $\sigma(\beta)$ are the same for all
$\alpha$ and $\beta$ in $C$. Let $\e_{\sigma}$ be the sign of
$\sigma(\alpha)$ for each real embedding $\sigma$. Then the
product over all real embeddings of $\e_\sigma$ is equal to $\e(C)$.
\qed

\section{Multiple Dedekind zeta functions}
\label{sec MDZF}
\subsection{Dedekind polylogarithms}
\label{subsec Dpolylog}

Let us recall the Dedekind zeta values
$$\zeta_K(m)=\sum_{\mathfrak{a}\neq(0)}\frac{1}{N_{K/\Q}(\mathfrak{a})^{m}},$$
where $\mathfrak{a}$ is an ideal in ${\cal{O}}_K$.

We are going to express the summation over elements, which belong to a finite union of
positive unimodular simple cones. We will define a Dedekind polylogarithm associates to a positive unimodular simple cone. The key result in this subsection will be that a Dedekind zeta value can be expressed as a $\Q$ linear combination of values of the Dedekind polylogarithms.

We also define a partial Dedekind zeta function
by summing over ideals in a given ideal class $[\mathfrak{a}]$
$$\zeta_{K,[\mathfrak{a}]}(m)=\sum_{\mathfrak{b}\in [\mathfrak{a}]}N_{K/\Q}(\mathfrak{b})^{-m},$$

Let us consider a partial Dedekind zeta functions
$\zeta_{K,[\mathfrak{a}]^{-1}}(m)$, corresponding to an ideal
class $[\mathfrak{a}]^{-1}$, where $\mathfrak{a}$ is an integral
ideal. For every integral ideal  $\mathfrak{b}$ in the class
$[\mathfrak{a}]^{-1}$, we have that
$$\mathfrak{a}\mathfrak{b}=(\alpha),$$
where $\alpha\in\mathfrak{a}$. Then
$$N_{K/\Q}(\mathfrak{b})=N_{K/\Q}(\mathfrak{a})^{-1}N_{K/\Q}((\alpha)).$$
Let
$$M(\mathfrak{a})=\bigcup_{i=1}^{n(\mathfrak{a})}C_i(\mathfrak{a}),$$
where $n(\mathfrak{a})$ is a positive integer and
$C_i(\mathfrak{a})$'s are unimodular simple cones.
Let $\alpha_{i}$ 
be an element of the intersection of
 ${\cal{O}}_K$ with  the dual cone 
 $C_i(\mathfrak{a})^*$, then 
 $\alpha_{i}C_i(\mathfrak{a})$ is a positive unimodular simple cone (see Lemma 
 \ref{lemma positive cone}).

Then,
\begin{align}
\nonumber \zeta_{K,[\mathfrak{a}]^{-1}}(m)
&=\sum_{\mathfrak{b}\in[\mathfrak{a}]^{-1}}N_{K/\Q}(\mathfrak{b})^{-m}=\\
\label{eq partial zeta}
&=N_{K/\Q}(\mathfrak{a})^m\sum_{i=1}^{n(\mathfrak{a})}
\e(C_i(\mathfrak{a}))^mN(\alpha_{i})^m\sum_{\a\in \alpha_{i}C_i(\mathfrak{a})}N_{K/\Q}(\alpha)^{-m},
\end{align}
where  $\e(C_i(\mathfrak{a}))=\pm 1$, depending on the cone,
$N_{K/\Q}(\mathfrak{a})$ is a norm of the ideal $\mathfrak{a}$ and 
$N(\alpha_{i})$ is the norm of the algebraic integer $\alpha_{i}$.

 We are going to give an example of higher
dimensional iteration in order to illustrate the usefulness of
this procedure. For a positive unimodular simple cone $C$, we define
$$f_{m}(C;u_1,\dots,u_n)=
\int^{u_1}_\infty\dots\int^{u_n}_\infty f_{m-1}(C;t_1,\dots,t_n) \di
t_1\wedge\dots\wedge \di t_n,$$ where $t_i\in(u_i,+\infty)$. This is an iteration, giving the simplest
type of iterated integrals on a membrane. We start the induction on $m$ from $m=0.$
Recall that $f_0$ was introduced in Definition \ref{def f0}.

Note that a norm of an algebraic number $\alpha$ can be expresses as a product of its embeddings in the complex numbers $\sigma_1(\alpha),\dots,\sigma_n(\alpha).$
$$N_{K/\Q}(\alpha)=\sigma_1(\alpha)\dots\sigma_n(\alpha).$$
Integrating term by term, we can express $f_m$ as an infinite sum
$$f_m(C;t_1,\dots,t_n)
=\sum_{\alpha\in C}\frac{\exp(-\sum_{i=1}^n \sigma_i(\alpha)t_i)}{N_{K/\Q}(\alpha)^m}.$$
Note that a cone $C$ is a linear combination of its generators so that the coefficients of the generators are positive integers. In particular, $0$ is not an element of an unimodular simple cone $C$, since then the generators are linearly independent over $\Q$. Thus, there is no division by $0$.

\begin{definition}
\label{def Dpolylog} We define an $m$-th  Dedekind polylogarithm,
associated to a number field $K$ and a positive unimodular simple cone
$C$, to be
$$Li^K_m(C;X_1,\dots,X_n)=f_m(C;-\log(X_1),\dots,-\log(X_n)).$$
\end{definition}
\begin{theorem}
\label{thm DZ} Dedekind zeta value at $s=m>1$ can be written as a
finite $\Q$-linear combination of Dedekind polylogarithms
evaluated at $(X_1,\dots,X_n)=(1,\dots,1)$.
\end{theorem}
\proof If $\mathfrak{a}_1,\dots,\mathfrak{a}_h$ are integral ideals
in ${\cal{O}}_K$, representing all the ideal classes, then using Equation \eqref{eq partial zeta}, we obtain
$$\zeta_K(m)=\sum_{j=1}^hN_{K/\Q}(\mathfrak{a}_j)^m
\sum_{i=1}^{n(\mathfrak{a})_j}\e(C_i(\mathfrak{a}_j))^mN(\alpha_{i,j})^m
f_m(\alpha_{i,j}C_i(\mathfrak{a}_j),0,\dots,0),
$$
where $C_i(\mathfrak{a}_j)$ are unimodular simple cones such that
$$\bigcup_{i=1}^{n(\mathfrak{a})_j}C_i(\mathfrak{a}_j)=M(\mathfrak{a_j})$$
and $\e(C_i(\mathfrak{a}_j))=\pm1$, depending on the cone 
(see Definition \ref{def fund domain} and Lemma \ref{lemma fund domain}). Let $\alpha_{i,j}\in C_i(\mathfrak{a}_j)^*\cap {\cal{O}}_K$ be an algebraic integer in the dual cone of 
$C_i(\mathfrak{a}_j)$. Then by Lemma \ref{lemma positive cone} we have that
 $\alpha_{i,j} C_i(\mathfrak{a}_j)$ is a positive unimodular simple cone.
The iterated integrals are hidden in the functions $f_m$. Consider Definition 
\ref{def it int n-forms}  with differential forms
$$\omega_1=f_0(\alpha_{i,j}C_i(\mathfrak{a}_j);z_1,\dots,z_n)\di z_1\wedge\dots\wedge \di z_n,$$
$$\omega_2=\omega_3=\dots=\omega_m=\di z_1\wedge\dots\wedge \di z_n.$$
And let $g$ be inclusion of $(0,\infty)^n$ in $\C^n$. Then the
corresponding iterated integral on a membrane gives
$$f_m(\alpha_{i,j}C_i(\mathfrak{a}_j);t_1,\dots,t_n).$$
\qed
\subsection{Multiple Dedekind zeta values}
\label{subsec MDZV}
We recall an integral representation of  a multiple zeta value
\[\zeta(k_1,k_2,\dots,k_m)=\sum_{0<n_1<\cdots<n_m}\frac{1}{n_1^{k_1}\cdots n_m^{k_m}}\]
has the following integral representation (see for example \cite{G}):
\begin{align}
\label{eq MZV}
\zeta(k_1,k_2,\dots,k_m)&=\int_{0<x_1<\cdots<x_{k_1+\cdots+k_m}}
\frac{dx_1}{1-x_1}\wedge\left(\frac{dx_2}{x_2}\wedge\cdots\wedge\frac{dx_{k_1}}{x_{k_1}}\right)
\wedge\\
\nonumber
&\wedge\frac{dx_{k_1+1}}{1-x_{k_1+1}}\wedge\left(\frac{dx_{k_1+2}}{x_{k_1+2}}\wedge\cdots\wedge\frac{dx_{k_1+k_2}}{x_{k_1+k_2}}\right)\wedge\\
\nonumber
&\cdots\\
\nonumber
&\wedge\frac{dx_{k_1+\cdots+k_{m-1}+1}}{1-x_{k_1+\cdots+k_{m-1}+1}}\wedge\left(\frac{dx_{k_1+\cdots+k_{m-1}+2}}{x_{k_1+\cdots+k_{m-1}+2}}\wedge\cdots\wedge\frac{dx_{k_1+\cdots+k_m}}{x_{k_1+\cdots+k_m}}\right)
\end{align}

Note that there are essentially two types of differential $1$-forms under the integral: 
$dx/(1-x)$ and $dx/x$. If we set $x_i=e^{-t_i}$, then we obtain the following formula, needed for the generalization to multiple Dedekind zeta values:

\begin{align}
\label{eq MZVt}
\zeta(k_1,k_2,\dots,k_m)&=\int_{t_1>\cdots>t_{k_1+\cdots+k_m}>0}
\frac{dt_1}{e^{t_1}-1}\wedge\left(dt_2\wedge\cdots\wedge dt_{k_1}\right)
\wedge\\
\nonumber
&\wedge\frac{dt_{k_1+1}}{e^{t_{k_1+1}}-1}\wedge\left(dt_{k_1+2}\wedge\cdots
\wedge dt_{k_1+k_2}\right)\wedge\\
\nonumber
&\cdots\\
\nonumber
&\wedge\frac{dt_{k_1+\cdots+k_{m-1}+1}}{e^{x_{k_1+\cdots+k_{m-1}+1}}-1}\wedge\left(dt_{k_1+
\cdots+k_{m-1}+2}\wedge\cdots\wedge dt_{k_1+\cdots+k_m}\right)
\end{align}

Note that in Equation \eqref{eq MZVt}, we have used $m$ copies of $dt/(e^t-1)$. To find their order, first we shuffle a set $S$ with $m$ elements (corresponding to $m$ copies of $dt/(e^t-1)$) with another set $S_1$ consisting of $m_1=-m+k_1+\cdots+k_m$ elements (corresponding to $m_1$ copies of $dt$). We choose a shuffle $\tau_1\in Sh(m,m_1)$, such that $\tau_1(1)=1$. The reason is that the first differential form in the iterated integral in Equation \ref{eq MZV} has to be $dx/(1-x)$, which is needed for convergence. The corresponding $1$ forms in Equation \eqref{eq MZVt} is $dt/(e^t-1)$ The relation between the shuffle $\tau_1$ and the set of integers $k_1,\dots,k_m$ is the following:
\begin{align*}
1&=\tau_1(1)\\
k_1+1&=\tau_1(2)\\
k_1+k_2+1&=\tau_1(3)\\
&\cdots\\
k_1+\cdots+k_{m-1}+1&=\tau_1(m)\\
k_1+\cdots+k_{m-1}+k_m&=\text{number of differential 1-forms}\\
\end{align*}
The integers $1,k_1+1,k_1+k_2+1,\dots,k_1+\cdots+k_{m-1}+1$, are the values of the index $i$, where the analogue of the form $dt_i/(e^{t_i}-1)$ appears under the integral, not the form $dt_i$.
(see Equation \eqref{eq MZVt})

In order to define multiple Dedekind zeta values, we will use $n$ shuffles of pairs of ordered sets, 
where $n=[K:\Q]$ is the degree on the number field.

Let $m_1,\dots,m_n$ be positive integers. (The positive integer
$m_i$ will denote the number of times the differential form $\di z_i$ occurs.) We define the following ordered sets:
$$S=\{1,2,\dots,m\},$$
$$S_i=\{m+1,\dots,m+m_i\},$$

\begin{definition}
Denote by $Sh^1(p,q)$ the subset of all shuffles $\tau\in Sh(p,q)$ of the two sets $\{1,\dots,p\}$ and $\{p+1,\dots,p+q\}$ such that $\tau(1)=1$
\end{definition}

For the definition of multiple Dedekind zeta values at the
positive integers, we use Definition \ref{def it int n-forms and 1-forms}, where we take the $n$-forms to be
$$\omega_j=f_0(C_j,z_1,\dots,z_n)\di z_1\wedge\dots\wedge \di z_n,$$
for $j=1,\dots, m$, where $C_1,\dots,C_m$ are positive unimodular simple cones,
and the $1$-forms to be $\di z_i$ on $\C^n$ occurring $m_i$ times for $i=1,\dots,n$.

\begin{definition}
\label{def MDZV} (Multiple Dedekind zeta values) For each $i=1,\dots,n$, 
let $\tau_i\in Sh^1(m,m_i)$. 
We define the integers
 $k_{i,j}$ and $m_i$ in terms of  the shuffle $\tau_i$ via the following relations
\begin{align}
1&=\tau_i(1)\\
\label{eq 1}k_{i,1}+1&=\tau_i(2)\\
\label{eq 2}k_{i,1}+k_{i,2}+1&=\tau_i(3)\\
\nonumber&\cdots\\
\label{eq m-1}k_{i,1}+\cdots+k_{i,m-1}+1&=\tau_i(m)\\
\label{eq m}
k_{i,1}+\cdots+k_{i,m-1}+k_{i,m}&=m+m_i
\end{align}
We define multiple Dedekind zeta
values at the positive integers by
$$
\zeta_{K;C_1,\dots C_m}\left(k_{1,1},\dots,k_{1,m};\dots;k_{n,1},\dots,k_{n,m}\right)
=
\int_{(g,\tau)}\omega_1\dots\omega_m(\di z_1)^{m_1}\dots(\di z_n)^{m_n}
$$
\end{definition}

\begin{theorem}
\label{thm MDZV}
For the general form of a multiple Dedekind zeta value, we
need: a number field $K$;  positive unimodular simple cones $C_1,\dots,C_m$  in
${\cal{O}}_K$; elements $\alpha_{j}\in C_j$
for $j=1,\dots,m$; complex embeddings of the elements
$\alpha_{i,j}=\sigma_i(\alpha_j)$. Then a multiple
Dedekind zeta value has the following representation as an infinite sum
\begin{align}
\nonumber 
&\zeta_{K;C_1,\dots C_m}\left(k_{1,1},\dots,k_{1,m};\dots;k_{n,1},\dots,k_{n,m}\right)
=\\
&=
\sum_{\alpha_1\in C_1}\cdots\sum_{\alpha_m\in C_m}
\prod_{i=1}^n \prod_{j=1}^m\left(\alpha_{i,1}+\cdots+\alpha_{i,j}\right)^{-k_{i,j}}.
\end{align}
 \end{theorem}
\proof There are $n$ different embedding $\sigma_1,\dots,\sigma_n$ of $K$ into $\C$.
Given $k_{i,1},\dots,k_{i,m}$ we find $m_i$ using Equation \eqref{eq m}.
Then we find $\tau_i$ by the values at $1,2,\dots,m$ obtained from Equations \eqref{eq 1}, \eqref{eq 2}, \eqref{eq m-1}. 

Now we use Definition \ref{def MDZV} of a multiple Dedekind zeta value 
in terms of iterated integrals on a membrane from Definition \ref{def it int n-forms and 1-forms}. 
We are going to follow closely Equation \eqref{eq MZVt}. 
The variable $t_{1,1}$ enters as a variable of the function $f_0(C_1;\cdots)$, 
the variables $t_{1,2},\dots,t_{1,k_{1,1}}$ appear 
as differential $1$-forms $dt_{1,2},\dots,dt_{1,k_{11}}$, since $\tau_1(2)=k_{1,1}+1$ 
(see Equation \eqref{eq 1}).
Recall that $\sigma_1:K\rightarrow \C$ is an embedding of $K$ into the complex numbers and $\alpha_{1,j}=\sigma_1(\alpha_j)$.
Thus integrating with respect to $t_{1,1},\dots,t_{1,k_{1,1}}$ 
gives us denominators $\alpha_{1,1}^{k_{1,1}}$, associated to each $\alpha_1\in\ C_1$
Then $t_{1,k_{1,1}+1}$ enters as a variable in the function $f_0(C_2;\cdots)$, since 
$\tau_1(2)=k_{1,1}+1$. 
Then the variables  $t_{1,k_{1,1}+2},\dots,t_{1,k_{1,1}+k_{1,2}}$ appear 
as differential $1$-forms $dt_{1,k_{1,1}+2},\dots,dt_{1,k_{1,1}+k_{1,2}}$, 
since $\tau_1(3)=k_{1,1}+k_{1,2}+1$ 
(see Equation \eqref{eq 2}). Thus integrating with respect to $t_{1,1},\dots,t_{1,k_{1,1}}$ 
gives us a denominators $\alpha_{1,1}^{k_{1,1}}(\alpha_{1,1}+\alpha_{1,2})^{k_{1,2}}$,
associated to each $\alpha_1\in C_1$ and each $\alpha_2\in C_2$.
Continuing this process to the variable $t_{1,m+m_1}$, we obtain a denominator
\[ \prod_{j=1}^m\left(\alpha_{1,1}+\cdots+\alpha_{1,j}\right)^{k_{1,j}},\]
associated to each $m$-tuple $(\alpha_1,\dots,\alpha_m)$, where $\alpha_j\in C_j$.
There are $n$ different embeddings $\sigma_1,\dots,\sigma_n$  of $K$ into $\C$, where $n=[K:\Q]$ is the degree of the number field. So far we have considered the contribution of the first embedding. The contribution of the first and the second embedding is obtained in essentially the same way as for the first embedding. It gives a denominator
\[ \prod_{j=1}^m\left(\alpha_{1,1}+\cdots+\alpha_{1,j}\right)^{k_{1,j}}\left(\alpha_{2,1}+\cdots+\alpha_{2,j}\right)^{k_{2,j}},\]
associated to each $m$-tuple $(\alpha_1,\dots,\alpha_m)$, where $\alpha_j\in C_j$.
Similarly, after integrating with respect to all the variables $t_{i,j}$ we obtain a denominator 
 \[\prod_{i=1}^n \prod_{j=1}^m\left(\alpha_{i,1}+\cdots+\alpha_{i,j}\right)^{k_{i,j}},\]
 associated to each $m$-tuple $(\alpha_1,\dots,\alpha_m)$, where $\alpha_j\in C_j$.
 Then the numerators are all equal to $1$ since the lower bound for the variables under the exponents in $f_0(C_j;\cdots)$ is $0$. Thus, the exponents become equal to $1$. 

The following examples of MDZV give analogues of Dedekind zeta function and of
(multiple) Eisenstein-Kronecker series.

{\bf{Examples:}}
1.  Let $C$ be a positive unimodular simple cone in the ring of integers ${\cal{O}}_K$ of a number filed $K$. In $m=1$ and if all values $k_{i,1}$ are equal to $k$, then
$$\zeta_{K;C}(k,\dots,k)=\sum_{\alpha\in C}\frac{1}{N_{K/\Q}(\alpha)^k}.$$
Note that the number $0$ does not belong to any  unimodular simple cone $C$.

2. Let $m=2$, and let 
$$k=k_{1,1}=\cdots=k_{n,1}$$
and
$$l=k_{1,2}=\cdots=k_{n,2}$$
be positive integers greater that $1$. 
Finally, let $C_1$ and $C_2$ be positive unimodular simple cones in the
ring of integers ${\cal{O}}_K$ of a number field $K$. Then the
corresponding multiple Dedekind zeta value can be written both as
a sum and as an integral:
 \begin{eqnarray}
 \label{DZ}
\zeta_{K;C_1,C_2}(k,\dots,k;l,\dots,l)=
\sum_{\alpha\in C_1,\beta\in C_2}\frac{1}{N_{K/\Q}(\alpha)^{k}N_{K/\Q}(\alpha+\beta)^{l}}
\end{eqnarray}

3. Let $K$ be an imaginary quadratic field. Let $C$ be a positive
unimodular simple cone in ${\cal{O}}_K$. We can represent the cone
$C$ as an $\N$-module: $C=\N\{\mu,\nu\}$, for $\mu,\nu\in
{\cal{O}}_K$. Put $z=\mu/\nu$. Consider then
$$\zeta_{K;C}(k,k)
=\sum_{\alpha\in C}\frac{1}{N(\alpha)^{k}}=|\nu|^{-2k}\sum_{a,b\in \N}\frac{1}{|az+b|^{2k}},$$
where the last sum is a portion of the $k$-th Eisenstein-Kronecker series. 
Such series could be found in \cite{W}
$$E_k(z)
=
\sum_{a,b\in\Z; (a,b)\neq(0,0)}\frac{1}{|az+b|^{2k}}.$$

4. With the notation of Example 3, we obtain an analogue of values of
multiple Eisenstein-Kronecker series
$$\zeta_{K;C,C}(k,k;l,l)=|\nu|^{-k-l}\sum_{a,b,c,d\in \N}\frac{1}{|az+b|^{2k}|(a+c)z+(b+d)|^{2l}}.$$
An alternative generalization was considered in \cite{G2}, Section 8.2.

\subsection{Multiple Dedekind zeta functions}
\label{subsec MDZF}

We will try to give some intuition behind the integral representation of the multiple zeta functions (see in \cite{G}).
After that we will generalize the construction to define the number field analogues - multiple Dedekind zeta functions.
In order to do that, we give two examples - one for $\zeta(3)$ and another for $\zeta(1,3)$.

We have
\begin{align*}
\zeta(3)&=\int_{0<x_1<x_2<x_3<1}\frac{\di x_1}{1-x_1}\wedge\frac{\di x_2}{x_2}\wedge\frac{\di x_3}{t_3}=\\
&=\int_{t_1>t_2>t_3>0}\frac{\di t_1\wedge \di t_2\wedge \di t_3}{e^{t_1}-1}=\\
&=\int_0^\infty\frac{t^2_1\di t_1}{\Gamma(3)(e^{t_1}-1)}.
\end{align*}
The first equality is due to Kontsevich. 
The second equality uses the change of variables $x_i=e^{-t_i}$. 
Both representations were examined in more details in Section \ref{sec examples}.  
The last equality uses the following equation
\begin{equation}
\label{eq dt}
\int_{b>t_1>t_2>\dots>t_n>a}\di t_1\wedge \di t_2\wedge \dots\wedge\di t_n
=\frac{(b-a)^{n}}{\Gamma(n+1)},
\end{equation}
whose proof we leave for the reader.

Similarly,
\begin{align*}
\zeta(1,3)
&=\int_{0<x_1<x_2<x_3<x_4<1}
\frac{\di x_1}{1-x_1}\wedge\frac{\di x_2}{1-x_2}\wedge\frac{\di x_3}{x_3}\wedge \frac{dx_4}{x_4}=\\
&=\int_{t_1>t_2>t_3>t_4>0}\frac{\di t_1\wedge \di t_2\wedge \di t_3\wedge\di t_4}{(e^{t_1}-1)(e^{t_2}-1)}=\\
&=\int_{t_1>t_2>0}\frac{\di t_1}{\Gamma(1)(e^{t_1}-1)}\wedge \frac{t_2^{3-1}\di t_2}{\Gamma(3)(e^{t_2}-1)}=\\
&=\int_{(0,\infty)^2}\frac{u_1^{1-1} u_2^{3-1}\di u_1\wedge \di u_2}{\Gamma(1)\Gamma(3)(e^{u_1+u_2}-1)(e^{u_2}-1)}.
\end{align*}
The first two equalities are of the same type as in the previous example. 
For the third equality we use Equation \eqref{eq dt}. For the last equality we use the change of variable
\begin{align*}
t_2&=u_2,\\
t_1&=u_1+u_2,
\end{align*}
where $u_1>0$ and $u_2>0$.
Following \cite{G}, we can interpolate the multiple zeta values by
$$\zeta(s_1,\dots,s_d)=\Gamma(s_1)^{-1}\dots\Gamma(s_d)^{-1}
\int_{(0,\infty)^d}\frac{u_1^{s_1-1}\dots u_d^{s_d-1}\di u_1\wedge\dots\wedge\di u_d}{(e^{u_1+\dots + u_d}-1)(e^{u_2+\dots + u_d}-1)\dots(e^{u_d}-1)}.
$$

If we denote by 
$$f_0(\N;t)=\sum_{a\in \N}e^{-at},$$
then
$$f_0(\N,t)=\frac{1}{e^t-1}$$
and
$$
\zeta(s_1,\dots,s_d)=\Gamma(s_1)^{-1}\dots\Gamma(s_d)^{-1}
\int_{(0,\infty)^d}\bigwedge_{j=1}^d f_0(\N;u_i+\dots+u_d)u_j^{s_j-1}\di u_j
$$

Let $n=[K:\Q]$ be the degree of the number field.
We recall Definition \ref{def f0} of $f_0$,
$$f_0(C;t_1,t_2,\dots,t_n)=\sum_{\alpha\in C}e^{-\sum_{i=1}^n \sigma_i(\alpha)t_i},$$
where $\sigma_i:K\rightarrow\C$ run through all embeddings of the field $K$ into the complex numbers.
Let $C_1$ and $C_2$ be two unimodular simple cones.
We want to raise an algebraic integer $\alpha_1$  to a complex power $s_1$ as a portion of the multiple Dedekind zeta function. We define
\[\alpha_1^{s_1}=e^{s_1\log(\alpha_1)}\]
for one element $\alpha\in C$ and $\alpha_1=\sigma_1(\alpha)$.
Choose a branch of the logarithmic function by making a cut of the complex plane at the negative real numbers. Since $C$ is a positive unimodular simple cone we have that $\R^n_{>}\subset C^*$, the function $\sigma_i$ composed with $\log$ is well defined on a positive unimodular simple cone $C$.

Then, we define a double Dedekind zeta function as
\begin{align}
&\zeta_{K;C_1,C_2}(s_{1,1},\dots,s_{n,1};s_{1,2},\dots,s_{n,2})
=
\\
\nonumber
&=\Gamma(s_{1,1})^{-1}\dots\Gamma(s_{n,2})^{-1}\int_{(0,+\infty)^{2n}}f_0(C_1;(u_{1,1}+u_{1,2}),\dots,(u_{n,1}+u_{n,2}))\times
\\
\nonumber
&\times f_0(C_2;(u_{1,2},\dots,u_{n,2})\bigwedge_{i=1}^n u_{i,1}^{s_{i,1-1}}\di u_{i,1}
\wedge
\bigwedge_{i=1}^n u_{i,2}^{s_{i,2}}\di u_{i,2}.
\end{align}
This definition combines both double zeta function and multiple Dedekind zeta values with double iteration. More generally, we can interpolate all multiple Dedekind zeta values into multiple Dedekind zeta functions so that multiple zeta functions are particular cases. Again, we define $u_{i,j}^{s_{i,j}}$ by
\[u_{i,j}^{s_{i,j}}=e^{s_{i,j}\log(u_{i,j})},\]
along the branch of logarithm described above. 
\begin{definition}
\label{def MDZF}
(Multiple Dedekind zeta functions)
Let $n=[K;\Q]$ be the degree of the number field.
Let $C_1,\dots,C_m$ be $m$ positive unimodular simple cones in ${\cal{O}}_K$.
Let $u_{i,j}\in(0,\infty)$ for $i=1,\dots,n$ and $j=1,\dots,m$.
We define multiple Dedekind zeta functions by the integral
\begin{align}
\nonumber 
&\zeta_{K;C_1,\dots,C_m}(s_{1,1},\dots,s_{n,1};\dots;s_{1,m},\dots,s_{n,m})
 =\\
&=\prod_{(i,j)=(1,1)}^{(n,m)}\Gamma(s_{i,j})^{-1}\times
\\
&\times\int_{(0,+\infty)^{mn}}
\bigwedge_{j=1}^{m}
f_0(C_j;(u_{1,j}+\dots+u_{1,m}),\dots,(u_{n,j}+\dots+u_{n,m}))
\bigwedge_{i=1}^{n}
u_{i,j}^{s_{i,j}-1}
\di u_{i,j},
\nonumber
\end{align}
when $Re(s_{i,j})>1$.
\end{definition}

\begin{theorem} (Infinite Sum Representation)
\label{thm MDZF}
For the general form of a multiple Dedekind zeta function, we
need: a number field $K$; positive unimodular simple cones $C_j$  in
${\cal{O}}_K$, for $j=1,\dots,m$; elements $\alpha_{j}\in C_j$
for $j=1,\dots,m$; complex embeddings of the elements
$\alpha_{i,j}=\sigma_i(\alpha_j);$ Then, a multiple
Dedekind zeta function has the following infinite sum representation
\begin{align*}
&\zeta_{K;C_1,\dots,C_m}(s_{1,1},\dots,s_{n,1};\dots;s_{1,m},\dots,s_{n,m})=\\
&=
\sum_{\alpha_1\in C_1}\cdots\sum_{\alpha_m\in C_m}
\prod_{i=1}^n \prod_{j=1}^m\left(\alpha_{i,1}+\dots+\alpha_{i,j}\right)^{-s_{i,j}},
\end{align*}
when  $Re(s_{i,j})>1$.
\end{theorem}
\proof 
We have
\begin{align*}
&\zeta_{K;C_1,\dots,C_m}(s_{1,1},\dots,s_{n,1};\dots;s_{1,m},\dots,s_{n,m})=\\
&=\prod_{(i,j)=(1,1)}^{(n,m)}\Gamma(s_{i,j})^{-1}\times
\\
&\times\int_{(0,+\infty)^{mn}}
\bigwedge_{j=1}^{m}
f_0(C_j;(u_{1,j}+\dots+u_{1,m}),\dots,(u_{n,j}+\dots+u_{n,m}))
\bigwedge_{i=1}^{n}
u_{i,j}^{s_{i,j}-1}
\di u_{i,j}=\\
&=\left(\prod_{(i,j)=(1,1)}^{(n,m)}\Gamma(s_{i,j})\right)^{-1}\sum_{\alpha_1\in C_1}\cdots\sum_{\alpha_m\in C_m}
\int_{(0,+\infty)^{mn}}
 \bigwedge_{j=1}^m\bigwedge_{i=1}^n e^{-(\alpha_{i,1}+\dots+\alpha_{i,j})u_{i,j}}u_{i,j}^{s_{i,j}-1}du_{i,j}\\
&=
\sum_{\alpha_1\in C_1}\cdots\sum_{\alpha_m\in C_m}
\prod_{i=1}^n \prod_{j=1}^m\left(\alpha_{i,1}+\dots+\alpha_{i,j}\right)^{-s_{i,j}}.
\end{align*}

The following examples give a bridge between Dedekind zeta function and values of Eisenstein series (Example 5), and between multiple Dedekind zeta function
and values of multiple Eisenstein series. More about values of multiple Eisenstein series will appear in Examples 7, 8, 9 on pages 29 and 30.

{\bf{Examples:}}

5. Let $K$ be any number field, let $m=1$ and let $C$ be 
a positive unimodular simple cone in ${\cal{O}}_K$. Then
\begin{equation}
\label{eq no it}
\zeta_{K;C}(s_{1,1},\dots,s_{n,1})=\sum_{\alpha\in C}\frac{1}{\prod_{i=1}^n \alpha_i^{s_{i,1}}},
\end{equation}
where $\alpha_i=\sigma_i(\alpha)$ is the $i$-th embedding in the complex numbers.
In particular, if all variables $s_{i,1}$, for $i=1,\dots,n$ have the same value $s$, then
\begin{equation}
\label{eq on it Dedekind}
\zeta_{K;C}(s,\dots,s)=\sum_{\alpha\in C}\frac{1}{N_{K/\Q}(\alpha)^s}.
\end{equation}

6. Now, let $m=2$. Then we have a double iteration. Let $K$ be any
number field. Let $C_1$ and $C_2$ be two positive unimodular simple cones.
Then
\begin{equation}
\label{eq double it}
\zeta_{K;C_1,C_2}\left(s_{1,1},\dots,s_{n,1};s_{1,2},\dots,s_{n,2}\right)
=
\sum_{\alpha\in C_1, \beta\in C_2}
\frac{1}{\prod_{i=1}^n \alpha_i^{s_{i,1}}(\alpha_i+\beta_i)^{s_{i,2}} }.
\end{equation}

In particular, if \[s_j=s_{1,j}=\dots=s_{n,j}\]
for $j=1,2$, then
\begin{equation}
\label{eq double Dedekind}
\zeta_{K;C_1,C_2}\left(s_1,\dots,s_1;s_2,\dots,s_2\right)
=
\sum_{\alpha\in C_1, \beta\in C_2}
\frac{1}{N_{K/\Q}(\alpha)^{s_{1}}N_{K/\Q}(\alpha+\beta)^{s_{2}} }.
\end{equation}


\section{Analytic properties and special values}
\label{sec analytic}

\subsection{Applications to multiple Eisenstein series}
\label{subsec Eis}
Assuming the analytic continuation (Theorem \ref{thm analytic}), we can consider values of the multiple Dedekind zeta functions, when one or more of the arguments are zero, which allows us to express special values of multiple Eisenstein series (see \cite{ZGK}) as multiple Dedekind zeta values. This is presented in the following three examples.

{\bf{Examples:}}
7. Let $K$ be an imaginary quadratic field. Let $C$ be a
positive unimodular simple cone in ${\cal{O}}_K$. We can represent the cone
$C$ as an $\N$-module: $C=\N\{\mu,\nu\}$, for $\mu,\nu\in
{\cal{O}}_K$. Put $z=\mu/\nu$. Consider
$\zeta_{K;C}(k_{1,1},k_{2,1})$ at $k_{1,1}=k$ and $k_{2,1}=0$.

Then
$$\zeta_{K;C}(k,0)
=\sum_{\alpha\in C}\frac{1}{\alpha_1^{k}}=\nu^{-k}\sum_{a,b\in \N}\frac{1}{(az+b)^k},$$
where the last sum is a portion of the $k$-th Eisenstein series.
$$E_k(\tau)
=
\sum_{a,b\in\Z; (a,b)\neq(0,0)}\frac{1}{(az+b)^k}.$$
is an analogue of Eisenstein series.

8. Let $K$ be an imaginary quadratic field. Let $C$ be a positive unimodular simple cone in 
${\cal{O}}_K$. We can represent $C$ as 
\[C=\N\{\mu,\nu\}=\{\alpha\in {\cal{O}}_K\mbox{ }|\mbox{ }\alpha=a\mu+b\nu,\mbox{ }a,b\in\N\}.\]
Put $z=\mu/\nu$. Then, we obtain a value of
multiple Eisenstein series
$$\zeta_{K;C,C}(k,0;l,0)=\nu^{-k-l}\sum_{a,b,c,d\in \N}\frac{1}{(az+b)^k((a+c)z+(b+d))^l}.$$

9. Similarly, one can define analogue of values of the above
Eisenstein series over real quadratic field $K$, by setting
$$E_k(z)=\nu^k\zeta_{K;C}(k,0)
=\sum_{\alpha\in C}\frac{1}{\alpha_1^{k}},$$
where $C=\N\{\mu,\nu\}$ is a positive unimodular simple cone in
a real quadratic ring of integers ${\cal{O}}_K$.

\subsection{Examples of Analytic continuation and Multiple Residues}
\label{subsec mult res}
The following examples of analytic continuations are based on a Theorem of Gelfand-Shilov. The constructions in example 11 is central for this Section. Using Example 11, we change the variables in a way that we can apply Gelfand-Shilov's Theorem (Theorem \ref{thm Gelfand-Shilov}) that gives analytic continuation of MDZF. Moreover, in Example 11 we compute a multiple residue at $(1,1,1,1)$. 
In Subsection \ref{sec remarks}, we generalize this method to other multiple residues and we state  two conjectures - one about the values of the multiple residues, again based on Examples  10 and 11, and other about more general MDZV.

Let us recall the theorem of Gelfand-Shilov.
 \begin{theorem}(\cite{GSh})
\label{thm Gelfand-Shilov}
Let $\phi(x)$ be a test function on $\R$, which decreases rapidly (exponentially) when $x\rightarrow \infty$ and let
\[x_{+}=\left\{
\begin{tabular}{cc}
$x$ & if $x>0$\\
$0$ & if  $x\leq 0$
\end{tabular}\right.\]
Then the value of the distribution $\frac{x_+^{s-1}dx}{\Gamma(s)}$ on the test function $\phi$, namely,
$$\int_\R\phi(x)\frac{x_+^{s-1}dx}{\Gamma(s)}$$
is an analytic function in the variable $s$.
\end{theorem}
In Examples 10 and 11, we express a multiple Dedekind zeta function (MDZF) as a test function times a distribution when $s_{i,j}>1$ up to $\Gamma$-factors. Then Theorem \ref{thm Gelfand-Shilov} tells us that we have an analytic continuation of the MDZF to all complex values of $s_{i,j}$ after multiplying by a suitable $\Gamma$-factors. Using this method, we compute the multiple residue at $(1,\dots,1)$.


{\bf{Example 10.}} 
Let $K$ be a quadratic field and let $C=\N\{\alpha,\beta\}$ be a positive unimodular simple cone. 
Put $\alpha_1$, $\alpha_2$ and $\beta_1$, $\beta_2$ 
be the images under the two embeddings into $\C$ of $\alpha$ and $\beta$, respectively.
We will compute the residue of 
$$\zeta_{K;C}(s_1,s_2)=\sum_{\mu\in C}\frac{1}{\mu_1^{s_1}\mu_2^{s_2}}$$
at the hyperplane $s_1+s_2=2$ and evaluated at $s_1=s_2=1$. 

We have
$$\zeta_{K;C}(s_1,s_2)=
\Gamma(s_1)^{-1}\Gamma(s_2)^{-1}
\int_0^\infty\int_0^\infty \frac{t_1^{s_1-1}t_2^{s_2-1}dt_1\wedge dt_2}
{(e^{\alpha_1t_1+\alpha_2t_2}-1)(e^{\beta_1t_1+\beta_2t_2}-1)}.$$
Set $t_1=x_1(1-x_2)$ and $t_2=x_1x_2.$
Then
\[\Gamma(s_1+s_2-2)^{-1}\zeta_{K;C}(s_1,s_2)\]
is the value of the distribution
$$Dx=\frac{x_{1+}^{s_1+s_2-3}x_{2+}^{s_2}(1-x_2)_+^{s_1-1}dx_1dx_2}
{\Gamma(s_1+s_2-2)\Gamma(s_1)\Gamma(s_2)}$$
at the test function
$$\phi=\frac{x_1^2}{(e^{x_1(\alpha_1+(\alpha_2-\alpha_1)x_2)}-1)(e^{x_1(\beta_1+(\beta_2-
\alpha_1)x_2)}-1)}.$$
Using Theorem \ref{thm Gelfand-Shilov} we obtain that $\Gamma(s_1+s_2-2)^{-1}\zeta_{K;C}(s_1,s_2)$ is an analytic function.
The residue of $\zeta_{K;C}(s_1,s_2)$ at $s_1+s_2=2$ is
$$\int_0^1\frac{1}{(\alpha_1+(\alpha_2-\alpha_1)x_2)(\alpha_1+(\beta_2-\beta_1)x_2)}\cdot
\frac{x_{2+}^{s_2}(1-x_2)_+^{s_1-1}dx_1dx_2}{\Gamma(s_1)\Gamma(s_2)}.$$
Evaluating at $s_2=1$, the integral becomes
$$\int_0^1\frac{dx_2}{(\alpha_1+(\alpha_2-\alpha_1)x_2)(\alpha_1+(\beta_2-\beta_1)x_2)}$$
Thus, the residue of $\zeta_{K;C}(s_1,s_2)$ at $s_1+s_2=2$, evaluated at $(s_1,s_2)=(1,1)$ is 
given by the above integral. After evaluating it, we obtain
$$\left.(Res_{s_1+s_2=2}\zeta_{K;C}(s_1,s_2))\right|_{(s_1,s_2)=(1,1)}=
\frac{\log\left(\frac{\alpha_2}{\alpha_1}\right)-\log\left(\frac{\beta_2}{\beta_1}\right)}
{\left|\begin{tabular}{ll}
$\alpha_1$ & $\beta_1$\\
$\alpha_2$ & $\beta_2$
\end{tabular}
\right|
}.$$
In particular, if $\beta=1$ we obtain
\begin{equation}
\label{eq res Dedekind}
\left.(Res_{s_1+s_2=2}\zeta_{K;C}(s_1,s_2))\right|_{(s_1,s_2)=(1,1)}=
\frac{\log(\alpha_2)-\log(\alpha_1)}{\alpha_2-\alpha_1}.
\end{equation}
Note that if $K$ is a real quadratic field and $\alpha$ is a generator of the group of units, then
$$|\log(\alpha_2)-\log(\alpha_1)|=2|\log(\alpha_1)|$$ 
is two times the regulator of the number field $K$ 
and $\alpha_2-\alpha_1$ is an integer multiple of the discriminant of $K$. For a definition of a discriminant and a regulator of a number field, one may consult \cite{IR}.
Equation \eqref{eq res Dedekind} is true for any quadratic field, not only for real quadratic fields.

The following Example gives key constructions needed for the proof of the analytic continuation of MDZF (Theorem \ref{thm analytic}). It is also a case study 
of Conjecture \ref{conj res} about the multiple residue of a multiple Dedekind zeta function at $(1,\dots,1)$.

{\bf{Example 11.}}
Let $K$ be a quadratic extension of $\Q$. Let $C_1=\N\{1,\alpha\}$ and $C_2=\N\{1,\gamma\}$ 
be two positive unimodular simple cones.
Let 
$$\zeta_{K;C_1,C_2}(s_1,s_2;s'_1,s'_2)
=
\sum_{\mu\in C_1;\mbox{ }\nu\in C_2}\frac{1}{\mu_1^{s_1}\mu_2^{s_2}(\mu_1+\nu_1)^{s'_1}(\mu_2+\nu_2)^{s'_2}}.$$
An integral representation can be written as
\begin{align}
\nonumber&\zeta_{K;C_1,C_2}(s_1,s_2;s'_1,s'_2)\Gamma(s_1)\Gamma(s_2)\Gamma(s'_1)\Gamma(s'_2)
=\\
\label{eq double Dedekind}&=
\int_{t_1>t'_1>0;\mbox{ }t_2>t'_2>0}
\frac{(t_1-t'_1)^{s_1-1}(t_2-t'_2)_2^{s_2-1}(t'_1)^{s'_1-1}(t'_2)^{s'_2-1}
dt_1\wedge dt_2\wedge dt'_1\wedge dt'_2 }
{(e^{t_1+t_2}-1)(e^{\alpha_1t_1+\alpha_2t_2}-1)(e^{t'_1+t'_2}-1)(e^{\gamma_1t'_1+\gamma_2t'_2}-1)}.
\end{align}
We compute the residue of a double Dedekind zeta function by taking the multiple residues of six functions $\zeta^{(a)},\dots,\zeta^{(f)}$ and considering their sum.

We shall write the last differential form in Equation \eqref{eq double Dedekind} as $Dt$.
Then, we have
\begin{equation}
\label{eq double Dedekind}
\zeta_{K;C_1,C_2}(s_1,s_2;s'_1,s'_2)
=(\Gamma(s_1)\Gamma(s_2)\Gamma(s'_1)\Gamma(s'_2))^{-1}
\int_{t_1>t'_1>0;\mbox{ }t_2>t'_2>0}
Dt
\end{equation}
We are going to express the above integral as a sum of six integrals, which are enumerated by all possible shuffles of $t_1>t'_1>0$ and $t_2>t'_2>0$. 
In other words, each of the six new integrals will be associated 
to each linear order among the variables $t_1,t_2,t'_1,t'_2$ 
that respect the above two inequalities among them.  
Thus, all possible cases are

(a) $t_1>t'_1>t_2>t'_2>0$,

(b) $t_1>t_2>t'_1>t'_2>0$,

(c) $t_1>t_2>t'_2>t'_1>0$,

(d) $t_2>t_1>t'_1>t'_2>0$,

(e) $t_2>t_1>t'_2>t'_1>0$,

(f) $t_2>t'_2>t_1>t'_1>0$.

Then the domain $t_1>t'_1>0;\mbox{ }t_2>t'_2>0$ can be represented as a disjoint union of the domains of integration, given in parts (a),...,(f). 
Then
\begin{align*}
\label{eq double Dedekind a-f}
\zeta_{K;C_1,C_2}(s_1,s_2;s'_1,s'_2)
&=(\Gamma(s_1)\Gamma(s_2)\Gamma(s'_1)\Gamma(s'_2))^{-1}
\int_{t_1>t'_1>0;\mbox{ }t_2>t'_2>0}
Dt\\
&=(\Gamma(s_1)\Gamma(s_2)\Gamma(s'_1)\Gamma(s'_2))^{-1}\times\\
&\times\left(
\int_{t_1>t'_1>t_2>t'_2>0}+
\int_{t_1>t_2>t'_1>t'_2>0}+\right.\\
&+\int_{t_1>t_2>t'_2>t'_1>0}+
\int_{t_2>t_1>t'_1>t'_2>0}+\\
&\left.+\int_{t_2>t_1>t'_2>t'_1>0}+
\int_{t_2>t'_2>t_1>t'_1>0}
\right)Dt=\\
&=
\zeta^{(a)}_{K;C_1,C_2}(s_1,s_2;s'_1,s'_2)+
\zeta^{(b)}_{K;C_1,C_2}(s_1,s_2;s'_1,s'_2)+\\
&+\zeta^{(c)}_{K;C_1,C_2}(s_1,s_2;s'_1,s'_2)+
\zeta^{(d)}_{K;C_1,C_2}(s_1,s_2;s'_1,s'_2)+\\
&+\zeta^{(e)}_{K;C_1,C_2}(s_1,s_2;s'_1,s'_2)+
\zeta^{(f)}_{K;C_1,C_2}(s_1,s_2;s'_1,s'_2).
\end{align*}
We define $\zeta^{(a)},\dots,\zeta^{(f)}$ to be the above six integrals, corresponding to the domains of integration given by (a),...,(f). The reason for defining them is to take multiple residues of the multiple Dedekind zeta function. It is easier to work with the functions 
$\zeta^{(a)},\dots,\zeta^{(f)}$ for the purpose of proving analytic continuation and taking residues. 

Thus, we compute the residue of a double Dedekind zeta function by taking the multiple residues of six functions $\zeta^{(a)},\dots,\zeta^{(f)}$ and considering their sum.

Consider  the domain of integration 
$$(a)\mbox{ } t_1>t'_1>t_2>t'_2>0.$$
We will compute the residues of 
\begin{align*}&\zeta_{K;C_1,C_2}^{(a)}(s_1,s_2;s'_1,s'_2)
\Gamma(s_1)\Gamma(s_2)\Gamma(s'_1)\Gamma(s'_2)
=\\
&=
\int_{t_1>t'_1>0>t_2>t'_2>0}
\frac{(t_1-t'_1)^{s_1-1}(t_2-t'_2)_2^{s_2-1}(t'_1)^{s'_1-1}(t'_2)^{s'_2-1}
dt_1\wedge dt_2\wedge dt'_1\wedge dt'_2 }
{(e^{t_1+t_2}-1)(e^{\alpha_1t_1+\alpha_2t_2}-1)(e^{t'_1+t'_2}-1)(e^{\gamma_1t'_1+\gamma_2t'_2}-1)}.
\end{align*}
Let the successive differences be $u_1=t_1-t'_1$, $u_2=t'_1-t_2$, $u_3=t_2-t'_2$, $u_4=t'_2$. Their admissive values are in the interval $(0,+\infty)$. Let us make the following substitution
\begin{align*}
u_1&=x_1(1-x_2),\\
u_2&=x_1x_2(1-x_3),\\
u_3&=x_1x_2x_3(1-x_4),\\
u_4&=x_1x_2x_3x_4.\\
\end{align*}
We are going to express the above integral in terms of the variables $x_1,$ $x_2,$ $x_3$ and $x_4$.
We have
\begin{align*}
\alpha_1t_1+\alpha_2t_2&=\alpha_1(u_1+u_2+u_3+u_4)+\alpha_2(u_3+u_4)
=x_1(\alpha_1+\alpha_2x_2x_3,)\\
t_1+t_2&=x_1(1+x_1x_3),\\
\gamma_1t'_1+\gamma_2t'_2&=\gamma_1(u_2+u_3+u_4)+\gamma_2u_4
=x_1x_2(1+x_3x_4),\\
t'_1+t'_2&=x_1x_2(1+x_3x_4),\\
t_1-t'_1&=u_1=x_1(1-x_2),\\
t_2-t'_2&=u_3=x_1x_2x_3(1-x_4),\\
t'_1&=u_2+u_3+u_4=x_1x_2,\\
t'_2&=u_4=x_1x_2x_3x_4.\\
\end{align*}
For the change of variables in the differential forms, we have
$$dt_1\wedge dt_2\wedge dt'_1\wedge dt'_2
=
du_1 \wedge du_2 \wedge du_3 \wedge du_4$$
and
$$\frac{du_1}{u_1}\wedge\frac{du_2}{u_2}\wedge\frac{du_3}{u_3}\wedge\frac{du_4}{u_4}
=
\frac{dx_1\wedge dx_2 \wedge dx_3 \wedge dx_4}{x_1x_2x_3x_4(1-x_2)(1-x_3)(1-x_4)}.$$
Then
\begin{align*}
&\zeta_{K;C_1,C_2}^{(a)}(s_1,s_2;s'_1,s'_2)
\Gamma(s_1)\Gamma(s_2)\Gamma(s'_1)\Gamma(s'_2)
=\\
&=
\int_{t_1>t'_1>0>t_2>t'_2>0}
\frac{(t_1-t'_1)^{s_1-1}(t_2-t'_2)_2^{s_2-1}(t'_1)^{s'_1-1}(t'_2)^{s'_2-1}
dt_1\wedge dt_2\wedge dt'_1\wedge dt'_2 }
{(e^{t_1+t_2}-1)(e^{\alpha_1t_1+\alpha_2t_2}-1)(e^{t'_1+t'_2}-1)(e^{\gamma_1t'_1+\gamma_2t'_2}-1)}\\
&=\int_{(u_1,\dots,u_4)\in(0,\infty)^4}
\frac{u_1^{s_1-1}u_3^{s_2-1}(u_2+u_3+u_4)^{s'_1-1}u_4^{s'_2-1}du_1 \wedge du_2 \wedge du_3 \wedge du_4}{(e^{t_1+t_2}-1)(e^{\alpha_1t_1+\alpha_2t_2}-1)(e^{t'_1+t'_2}-1)(e^{\gamma_1t'_1+\gamma_2t'_2}-1)}\\
&=\int_{(u_1,\dots,u_4)\in(0,\infty)^4}
\frac{u_1^{s_1}u_3^{s_2}(u_2+u_3+u_4)^{s'_1-1}u_4^{s'_2}u_2\frac{du_1}{u_1} \wedge \frac{du_2}{u_2} \wedge \frac{du_3}{u_3} \wedge \frac{du_4}{u_4}}{(e^{t_1+t_2}-1)(e^{\alpha_1t_1+\alpha_2t_2}-1)(e^{t'_1+t'_2}-1)(e^{\gamma_1t'_1+\gamma_2t'_2}-1)}\\
&=\int_0^\infty\int_{(0,1)^3}
\frac{[x_1(1-x_2)]^{s_1}[x_1x_2x_3(1-x_4)]^{s_2}[x_1x_2]^{s'_1}[x_1x_2x_3x_4]^{s'_2}[x_1x_2(1-x_3)]\Omega}{(e^{x_1(1+x_2x_3)-1})(e^{x_1(\alpha_1+\alpha_2x_2x_3)-1})(e^{x_1x_2(1+x_3x_4)}-1)(e^{x_1x_2(\gamma_1+\gamma_2x_3x_4)}-1)},
\end{align*}
where
$$\Omega
=\frac{du_1}{u_1} \wedge \frac{du_2}{u_2} \wedge \frac{du_3}{u_3} \wedge \frac{du_4}{u_4}
=\frac{dx_1\wedge dx_2 \wedge dx_3 \wedge dx_4}{x_1x_2x_3x_4(1-x_2)(1-x_3)(1-x_4)}.
$$
We can express the last integral as a distribution evaluated at a test function. Put
\begin{align*}
&\phi(x_1,x_2,x_3,x_4)=\\
&=\frac{x_1^4x_2^2}
{(e^{x_1(1+x_2x_3)-1})(e^{x_1(\alpha_1+\alpha_2x_2x_3)-1})(e^{x_1x_2(1+x_3x_4)}-1)(e^{x_1x_2(\gamma_1+\gamma_2x_3x_4)}-1)}
\end{align*}
to be a test function. Let $Dx$ be a distribution defined by
\begin{align} 
Dx=&
\frac{x_{1+}^{s_1+s_2+s'_1+s'_2-5}}{\Gamma(s_1+s_2+s'_1+s'_2-4)}
\frac{x_{2+}^{s_2+s'_1+s'_2-3}}{\Gamma(s_2+s'_1+s'_2-2)}\times\\
 &\times\frac{x_{3+}^{s_2+s'_2-1}}{\Gamma(s_2+s'_2)}
\frac{x_{4+}^{s'_2-1}}{\Gamma(s'_2)}
(1-x_2)_+^{s_1-1}
(1-x_4)_+^{s_2-1}.
\label{eq dist a}
\end{align}
Then we have
\begin{align*}
\zeta_{K;C_1,C_2}^{(a)}(s_1,s_2;s'_1,s'_2)=&\frac{\Gamma(s_1+s_2+s'_1+s'_2-4)\Gamma(s_2+s'_1+s'_2-2)\Gamma(s_2+s'_2)}
{\Gamma(s_1)\Gamma(s_2)\Gamma(s'_1)}\times\\
&\times\int\phi Dx
\end{align*}
Using Theorem \ref{thm Gelfand-Shilov} we prove 
 analytic continuation of $\zeta_{K;C_1,C_2}^{(a)}(s_1,s_2;s'_1,s'_2)$ everywhere except at the poles of $\Gamma(s_1+s_2+s'_1+s'_2-4)\Gamma(s_2+s'_1+s'_2-2)\Gamma(s_2+s'_2)$.
 
We will take the residue at $s_1+s_2+s'_1+s'_2=4$ and then evaluate at 
$(s_1,s_2,s'_1,s'_2)=(1,1,1,1).$
Note that 
$$\phi(0,x_2,x_3,x_4)=
\frac{1}{(1+x_2x_3)(\alpha_1+\alpha_2x_2x_3)(1+x_3x_4)(\gamma_1+\gamma_2x_3x_4)}$$
Then we can compute
\begin{align*}
&(Res_{s_1+s_2+s_1'+s_2'=4}
\zeta^{(a)}_{K;C_1,C_2}(s_1,s_2;s_1',s_2'))|_{(1,1,1,1)}
=\\
&=
\int_{(0,1)^3}\frac{x_3dx_2\wedge dx_3\wedge dx_4} {(1+x_2x_3)(\alpha_1+\alpha_2x_2x_3)(1+x_3x_4)(\gamma_1+\gamma_2x_3x_4)}.
\end{align*}
Note that we cannot take a double residue at the point $(1,1,1,1)$. Also, the value that we obtain is a period.

Consider the domain in integration 
$$\mbox{(b) }t_1>t_2>t_1'>t_2'>0$$
We will compute the residues of
\begin{align*}
&\zeta^{(b)}_{K;C_1,C_2}(s_1,s_2;s_1',s_2')
\Gamma(s_1)\Gamma(s_2)\Gamma(s_1')\Gamma(s_2')
=\\
&=\int_{ t_1>t_2>t_1'>t_2'>0}
\frac{(t_1-t_1')^{s_1-1}(t_2-t_2')^{s_2-1}(t_1')^{s_1'-1}(t_2')^{s_2'-1}dt_1\wedge dt_2\wedge dt_1'\wedge dt_2'}
{(e^{t_1+t_2}-1)(e^{\alpha_1t_1+\alpha_2t_2}-1)(e^{t_1'+t_2'}-1)(e^{\gamma_1t_1'+\gamma_2t_2'}-1)}.
\end{align*}
Let the successive differences be
\begin{align*}
u_1&=t_1-t_2,\\
u_2&=t_2-t_1',\\
u_3&=t_1'-t_2',\\
u_4&=t_2'.
\end{align*}
Their admissive values are $(0,+\infty)$.
Let 
\begin{align*}
u_1&=x_1(1-x_2),\\
u_2&=x_1x_2(1-x_3),\\
u_3&=x_1x_2x_3(1-x_4),\\
u_4&=x_1x_2x_3x_4.
\end{align*}
We will express the above integral in terms of $x_1,\dots,x_4$.
We have
\begin{align*}
\alpha_1t_1+\alpha_2t_2&=\alpha_1(u_1+u_2+u_3+u+4)+\alpha_2(u_2+u_3+u_4)
=x_1(\alpha_1+\alpha_2x_2),\\
t_1+t_2&=x_1(1+x_2),\\
\gamma_1t_1'+\gamma_2t_2'&=\gamma_1(u_3+u_4)+\gamma_2u_4=x_1x_2x_3(\gamma_1+\gamma_2x_4),\\
t_1'+t_2'&=x_1x_2x_3(1+x_4),\\
t_1-t_1'&=u_1+u_2=x_1(1-x_2x_3),\\
t_2-t_2'&=u_2+u_3=x_1x_2(1-x_3x_4),\\
t_1'&=u_3+u_4=x_1x_2u_3,\\
t_2'&=u_4=x_1x_2x_3x_4.
\end{align*}
For  the differential forms,  we have
\begin{align*}
dt_1\wedge dt_2\wedge dt_1'\wedge dt_2'&=du_1\wedge du_2\wedge du_3\wedge du_4\\
\frac{du_1}{u_1}\wedge
\frac{du_2}{u_2}\wedge
\frac{du_3}{u_3}\wedge
\frac{du_4}{u_4}
&=
\frac{dx_1\wedge dx_2\wedge dx_3\wedge dx_4}{x_1x_2x_3x_4(1-x_2)(1-x_3)(1-x_4)}.
\end{align*}
Then
\begin{align*}
&\zeta^{(b)}_{K;C_1,C_2}(s_1,s_2;s_1',s_2')
\Gamma(s_1)\Gamma(s_2)\Gamma(s_1')\Gamma(s_2')
=\\
&=\int_{ t_1>t_2>t_1'>t_2'>0}
\frac{(t_1-t_1')^{s_1-1}(t_2-t_2')^{s_2-1}(t_1')^{s_1'-1}(t_2')^{s_2'-1}dt_1\wedge dt_2\wedge dt_1'\wedge dt_2'}
{(e^{t_1+t_2}-1)(e^{\alpha_1t_1+\alpha_2t_2}-1)(e^{t_1'+t_2'}-1)(e^{\gamma_1t_1'+\gamma_2t_2'}-1)}\\
&=\int_{(u_1,\dots,u_4)\in(0,\infty)^4}
\frac{(u_1+u_2)^{s_1-1}(u_2+u_3)^{s_2-1}(u_3+u_4)^{s_1'-1}u_4^{s_2'-1}du_1\wedge du_2\wedge du_3\wedge du_4}
{(e^{t_1+t_2}-1)(e^{\alpha_1t_1+\alpha_2t_2}-1)(e^{t_1'+t_2'}-1)(e^{\gamma_1t_1'+\gamma_2t_2'}-1)}\\
&=\int_{(u_1,\dots,u_4)\in(0,\infty)^4}
\frac{(u_1+u_2)^{s_1-1}(u_2+u_3)^{s_2-1}(u_3+u_4)^{s_1'-1}u_4^{s_2'}
u_1u_2u_3\frac{du_1\wedge du_2\wedge du_3\wedge du_4}{u_1u_2u_3u_4}}
{(e^{t_1+t_2}-1)(e^{\alpha_1t_1+\alpha_2t_2}-1)(e^{t_1'+t_2'}-1)(e^{\gamma_1t_1'+\gamma_2t_2'}-1)}\\
&=\int_0^\infty\int_{(0,1)^3}\frac{[x_1(1-x_2x_3)]^{s_1-1}[x_1x_2(1-x_3x_4)]^{s_2-1}[x_1x_2x_3]^{s_1'-1}[x_1x_2x_3x_4]^{s_2'}X\Omega}
{(e^{x_1(1+x_2x_3)}-1)(e^{x_1(\alpha_1+\alpha_2x_2x_3)}-1)(e^{x_1x_2(1+x_3x_4)}-1)(e^{x_1x_2(\gamma_1+ \gamma_2x_3x_4)}-1)},
\end{align*}
where $$X=x_1^3x_2^2x_3(1-x_2)(1-x_3)(1-x_4)$$ and
$$\Omega=\frac{du_1\wedge du_2\wedge du_3\wedge du_4}{u_1u_2u_3u_4}
=
\frac{dx_1\wedge dx_2\wedge dx_3\wedge dx_4}{x_1x_2x_3x_4(1-x_2)(1-x_3)(1-x_4)}$$
Now we can express the last integral as a distribution evaluated at a test function.
Put 
\begin{align*}
&\phi(x_1,x_2,x_3,x_4)=\\
&\frac{x_1^4x_2^2x_3^2}
{(e^{x_1(1+x_2)}-1)(e^{x_1(\alpha_1+\alpha_2x_2)}-1)(e^{x_1x_2x_3(1+x_4)}-1)(e^{x_1x_2x_3(\gamma_1+ \gamma_2x_4)}-1)},
\end{align*}
to be a test function.
Let $Dx$ be a distribution defined by
\begin{align}
Dx&=\frac{x_{1+}^{s_1+s_2+s_1'+s_2'-5}}{\Gamma(s_1+s_2+s_1'+s_2'-4)}
\frac{x_{2+}^{s_2+s_1'+s_2'-3}}{\Gamma(s_2+s_1'+s_2'-2)}\times\\
\label{eq dist b}
&
\times
\frac{x_{3+}^{s_2+s_2'-3}}{\Gamma(s_2+s_2'-2)}
\frac{x_{4+}^{s_2'-1}}{\Gamma(s_2')}(1-x_2x_3)_+^{s_1-1}(1-x_3x_4)_+^{s_2-1}
\end{align}
Then we have
\begin{align*}
\zeta^{(b)}_{K;C_1,C_2}(s_1,s_2;s_1',s_2')=&
\frac{\Gamma(s_1+s_2+s_1'+s_2'-4)\Gamma(s_2+s_1'+s_2'-2)\Gamma(s_2+s_2'-2)}
{\Gamma(s_1)\Gamma(s_2)\Gamma(s_1')}\times\\
&\times\int\phi Dx
\end{align*}
Using Theorem \ref{thm Gelfand-Shilov} we prove 
 analytic continuation of $\zeta_{K;C_1,C_2}^{(b)}(s_1,s_2;s'_1,s'_2)$ everywhere except at the poles of $\Gamma(s_1+s_2+s_1'+s_2'-4)\Gamma(s_2+s_1'+s_2'-2)\Gamma(s_2+s_2'-2)$.

We will take the residues at $s_1+s_2+s_1'+s_2'=4$ and at $s_1'+s_2'=2$. And then, we will evaluate at 
$(s_1,s_2,s_1',s_2')=(1,1,1,1)$. 
Note that 
$$\phi(0,x_2,0,x_4)=
\frac{1}{(1+x_2)(\alpha_1+\alpha_2x_2)(1+x_4)(\gamma_1+\gamma_2x_4)}$$
Then we can compute
\begin{align*}
&(Res_{s_1'+s_2'=2}Res_{s_1+s_2+s_1'+s_2'=4}
\zeta^{(b)}_{K;C_1,C_2}(s_1,s_2;s_1',s_2'))|_{(1,1,1,1)}
=\\
&=
\int_{(0,1)^2}\frac{dx_2\wedge dx_4} {(1+x_2)(\alpha_1+\alpha_2x_2)(1+x_4)(\gamma_1+\gamma_2x_4)}=\\
&=
\frac{\log\left(\frac{2\alpha_1}{\alpha_1+\alpha_2}\right)}{\alpha_1-\alpha_2}
\cdot
\frac{\log\left(\frac{2\gamma_1}{\gamma_1+\gamma_2}\right)}{\gamma_1-\gamma_2}.
\end{align*}

For the cases (c), (d) and (e), we obtain
\begin{align*}
&(Res_{s_1'+s_2'=2}Res_{s_1+s_2+s_1'+s_2'=4}
\zeta^{(c)}_{K;C_1,C_2}(s_1,s_2;s_1',s_2'))|_{(1,1,1,1)}
=\\
&=
\int_{(0,1)^2}\frac{dx_2\wedge dx_4} {(1+x_2)(\alpha_1+\alpha_2x_2)(x_4+1)(\gamma_1x_4+\gamma_2)}=\\
&=
-\frac{\log\left(\frac{2\alpha_1}{\alpha_1+\alpha_2}\right)}{\alpha_1-\alpha_2}
\cdot
\frac{\log\left(\frac{2\gamma_2}{\gamma_1+\gamma_2}\right)}{\gamma_1-\gamma_2}.
\end{align*}

\begin{align*}
&(Res_{s_1'+s_2'=2}Res_{s_1+s_2+s_1'+s_2'=4}
\zeta^{(d)}_{K;C_1,C_2}(s_1,s_2;s_1',s_2'))|_{(1,1,1,1)}
=\\
&=
\int_{(0,1)^2}\frac{dx_2\wedge dx_4} {(x_2+1)(\alpha_1x_2+\alpha_2)(1+x_4)(\gamma_1+\gamma_2x_4)}=\\
&=
-\frac{\log\left(\frac{2\alpha_2}{\alpha_1+\alpha_2}\right)}{\alpha_1-\alpha_2}
\cdot
\frac{\log\left(\frac{2\gamma_1}{\gamma_1+\gamma_2}\right)}{\gamma_1-\gamma_2}.
\end{align*}

\begin{align*}
&(Res_{s_1'+s_2'=2}Res_{s_1+s_2+s_1'+s_2'=4}
\zeta^{(e)}_{K;C_1,C_2}(s_1,s_2;s_1',s_2'))|_{(1,1,1,1)}
=\\
&=
\frac{\log\left(\frac{2\alpha_2}{\alpha_1+\alpha_2}\right)}{\alpha_1-\alpha_2}
\cdot
\frac{\log\left(\frac{2\gamma_2}{\gamma_1+\gamma_2}\right)}{\gamma_1-\gamma_2}.
\end{align*}

Case (f) is similar to case (a), namely, there is no double residue at the point $(1,1,1,1)$.
Thus, we obtain
\begin{align*}
&Res_{s_1'+s_2'=2}Res_{s_1+s_2+s_1'+s_2'=4}
\zeta^{}_{K;C_1,C_2}(s_1,s_2;s_1',s_2'))|_{(1,1,1,1)}=\\
&=Res_{s_1'+s_2'=2}Res_{s_1+s_2+s_1'+s_2'=4}
(\zeta^{(b)}_{K;C_1,C_2}(s_1,s_2;s_1',s_2')+\zeta^{(c)}_{K;C_1,C_2}(s_1,s_2;s_1',s_2')+\\
&+\zeta^{(d)}_{K;C_1,C_2}(s_1,s_2;s_1',s_2')+\zeta^{(e)}_{K;C_1,C_2}(s_1,s_2;s_1',s_2'))|_{(1,1,1,1)}
=\\
&=
\frac{\log(\alpha_2)-\log(\alpha_1)}{\alpha_2-\alpha_1}
\cdot
\frac{\log(\gamma_2)-\log(\gamma_1)}{\gamma_2-\gamma_1}.
\end{align*}
Note that if $K$ is a real quadratic field and $\alpha$ is a generator of the group of units, then
$$|\log(\alpha_2)-\log(\alpha_1)|=2\log|\alpha_1|$$ 
is two times the regulator of the number field $K$ and
$$\alpha_2-\alpha_1$$ is an integer multiple of the discriminant of the field $K$. 
For a definition of a discriminant and a regulator of a number field one may consult with \cite{IR}. 
The above formula is true for any quadratic field, not necessarily for a real quadratic field.

\subsection{Analytic Continuation of Multiple Dedekind Zeta Functions}

\begin{theorem}
\label{thm analytic}
Multiple Dedekind zeta functions 
\[\zeta_{K;c_1,\dots,C_m}(s_{1,1},\dots,s_{n,1};\dots;s_{1,m},\dots,s_{n,m})\]
 have an analytic continuation from the region $Re(s_{i,j})>1$ for all $i$  and $j$ to $s_{i,j}\in \C$ with exception of hyperplanes. The hyperplanes are defined by  sum several of the variables $s_{i,j}$ without repetitions being set equal to an integer.
\end{theorem}
\proof
Recall that $f_0(C_j;t_{1,j},\dots,t_{n,j})$ is used to define the multiple Dedekind zeta functions, where the domain of
integration is
\[D=\{(t_{i,j})\in \R^{mn}|t_{i,1}>t_{i,2}>\dots >t_{n,j}>0\}.\]
\[J_i=\{k_{i,1},\dots,k_{i,m}\}.\]
Note that there are $n$ sets $J_1,\dots,J_n$, and each of them has $m$ elements,  $|J_i|=m$.
Let $\tau$ run through all the shuffles of the ordered sets $J_1,\dots,J_n$,
$$\tau\in Sh(J_1,\dots,J_n).$$
Let $t_1,t_2\dots,t_{mn}$ be the variables $t_{1,1},\dots,t_{n,m}$, written in decreasing order.  
There are finitely many ways of arranging the variables in decreasing order. More precisely, the number of such arrangements is equal to the number of shuffles in $Sh(J_1,\dots,J_n)$. We need to consider all such shuffles in order to express the multiple Dedekind zeta function as a sum of partial multiple Dedekind zeta functions, corresponding to each shuffle $\tau$ (see Example 11).
Let $u_k=t_k-t_{k+1}$, for $k=1,\dots,mn-1$ and 
$u_{mn}=t_{mn}.$
Let 
\begin{align*}
u_1&=x_1(1-x_2)\\
u_2&=x_1x_2(1-x_3)\\
&\cdots\\
u_{mn-1}&=x_1\dots x_{mn-1} (1-x_{mn})\\
u_{mn}&=x_1\dots x_{mn}.
\end{align*}
Note that \[(x_1,\dots,x_{mn})\in \R\times [0,1]^{mn-1}\]
Then each of the linear factors in the denominator of $f_0(C_j;\cdots)$ can be written as
$$x_1\dots x_kg_{j,l,\tau}.$$
for some positive integer $k$, $k\leq mn$ and a polynomial $g_{j,l,\tau}$,
 in the variables $x_1,\dots, x_{mn}$, vanishing at the origin. Also, the indices $i$ and $j$ of the polynomial $g_{j,l,\tau}$ are associated to the $l$-th generator of the cone $C_j$, and $\tau$ is a shuffle of  ordered sets $Sh(J_1,\dots,J_n)$. 
Example of polynomials $g_{j,l,\tau}$, 
can be found in Equations  \eqref{eq dist a} and \eqref{eq dist b}. 
If the polynomial $g_{j,l,\tau}$, in terms of $x_1,\dots,x_{mn}$, has a pole on the domain of integration $\R\times [0,1]^{mn-1}$
then we move the membrane of integration to a membrane $D'\in \C^{mn}$ so that the real coordinates of $D'$ give the domain $\R\times [0,1]^{mn-1}$. (This is analogue of contour integration.) The integral representation in terms of $x_1,\dots,x_{mn}$ (similar to the ones in Example 11, giving $\zeta^{(a)}$), are a type of zeta function that we call {\it{partial MDZF}} times $\Gamma$-factors. Using Theorem \ref{thm Gelfand-Shilov},
we find that the partial MDZF together with the $\Gamma$-factors is an analytic function. 
The $\Gamma$-factors give hyperplanes where the poles of the partial MDZF occur.
Expressing a MDZF as a finite sum of partial MDZF we obtain the analytic continuation from the domain $Re (s_{i,j})>1$ to $s_{i,j}\in \C$ with poles along hyperplanes coming from $\Gamma$-factors.
\qed

\subsection{Final remarks}
\label{sec remarks}
In this final Subsection, we proof that certain multiple residue of a  multiple Dedekind zeta functions is a period in the sense of algebraic geometry. Based on Theorem \ref{thm periods}, we state two conjectures. One of the conjectures is about the exact values of the multiple residue and the other conjecture is about values of the multiple Dedekind zeta functions at other integers.

\begin{theorem}
\label{thm periods}
The multiple residue of a multiple Dedekind zeta function at the point
\[(s_{1,1},\dots,s_{n,1};\dots;s_{1,m},\dots,s_{n,m})=(1,\dots,1)\]
is a period over $\Q$.
\end{theorem}
\proof
We use the notation introduced in the proof of Theorem \ref{thm analytic} and of Example 11 in Subsection \ref{subsec mult res}.

The $m$-fold residue at $(1,\dots,1)$ can be computed via an integral of a rational function, which is a product of the functions representing the hyperplanes, where $f_0$ vanishes, expressed in terms of the variables $x_i$. The value is a period. In general, after we take the multiple residues at $(1,\dots,1)$, we obtain an integral of a rational function, (which is a product of $g_{j,l,\tau}$ over $j$ and $l$, for $j=1,\dots,m$, where $l$ signifies the $l$-generator of the cone $C_j$). Note that $\tau$ is a shuffle. So that different shuffles $\tau$ correspond to different partial MDZF. The boundaries of the integral (after taking the multiple residues) form a unit cube. Therefore, the value of the multiple residue at $(1,\dots,1)$ of a multiple Dedekind zeta function is a period.
\qed

For a more precise interpretation see Conjecture \ref{conj res} and Examples 10 and 11.

From Examples 10 and 11, we know that a multiple residue of multiple Dedekind zeta function
is a product of residues of partial Dedekind zeta functions, for quadratic fields and double iteration. 
For unimodular simple cones $C_1,\dots,C_m$, we consider a multiple Dedekind zeta
function
$$\zeta_{K;C_1,\dots,C_m}(s_1,\dots,s_d)=\sum_{\alpha_1\in C_1,\dots\alpha_d\in C_m}
\frac{1}{N(\alpha_1)^{s_1}N(\alpha_1+\alpha_2)^{s_2}\cdots N(\alpha_1+\dots+\alpha_m)^{s_m}}.$$
We expect that
\begin{conjecture}
\label{conj res}
The multiple residue of $\zeta_{K;C_1,\dots,C_m}(s_{1},\dots,s_{m})$ at the point $(1,\dots,1)$, namely
$$Res_{s_m=1}\dots Res_{s_1+\dots s_m=m} \zeta_{K;C_1,\dots,C_m}(s_{1},\dots,s_{m})
=\prod_{j=1}^m Res_{s=1}\zeta_{K;C_j}(s).$$
\end{conjecture}
The conjecture is proven for  a quadratic fields $K$ and double iteration in Examples 10 and 11.

We do expect that  multiple
Dedekind zeta values should be periods over $\Q$.

\begin{conjecture}\label{conj periods}
Let $K$ be a number field. For any choice of unimodular simple cones $C_1,\dots,C_m$,
in the ring of integers of a number field $K$, we have that the multiple Dedekind zeta values
(see Definition \ref{def MDZV})
$$\zeta_{K;C_1,\dots,C_m}(k_{1,1},\dots,k_{n,1};\dots;k_{1,m},\dots,k_{n,m})$$
are periods over $\Q$ , when the  $k_{1,1},\dots,k_{m,n}$ are natural numbers greater than $1$.
\end{conjecture}

The reasons for this conjecture are the following: 
\begin{enumerate}
\item
We have that multiple zeta values are periods; 
\item 
Dedekind zeta values are periods; 
\item  
From Theorem \ref{thm periods}, we have that the multiple residue of a multiple Dedekind zeta function at $(1,\dots,1)$ is a period.

\item The main reason is the representation of multiple Dedekind zeta values as iterated integrals on membranes. We will give a semi-algebraic relations among the variables in such integrals. 
\end{enumerate} 
We use Equations \eqref{eq yj} and \eqref{eq f0}. Recall that 
\[y_j=\prod_{i=1}^n\exp(-\sigma_i(e_j)z_i).
\]
If we set 
\[x_i=\exp(-z_i),\]
then a multiple Dedekind zeta value is an iterated integral on a membrane of the $n$-form
\[\prod_j \frac{dy_j}{1-y_j}\]
and the $1$-forms 
\[\frac{dx_i}{x_i},\]
which mostly resembles polylogarithms. However, the relations between $x_i$'s and $y_j$'s are (semi)-algebraic, namely,
\[\frac{dy_j}{y_j}=\sum_{i=1}^n\sigma_i(e_j)\frac{dx_i}{x_i},\]
which are not algebraic. Explicitly, they are given by
\[\log(y_j)=\sum_{i=1}^n\sigma_i(e_j)\log(x_i),\]
which involves  the logarithmic function. Note that  the logarithmic function is a homotopy invariant function on a path space. In this setting, the above  logarithmic functions can be considered as a function on the path space of an affine $n$-space without some divisors. One may take  a simplicial scheme as a model of the path space so that it restricts well onto the loop space of a scheme as a simplicial scheme. Hopefully, that would interpret the above (semi)-algebraic relations in terms of logarithms in an algebraic context. 

{\bf{Acknowledgements:}} This paper owes a lot to many people.
Ronald Brown gave inspiring talk on higher cubical categories. I
few days after that talk, I had a good definition of iterated
integrals on a membrane. With Alexander Goncharov I had a lot of
fruitful discussions. His great interest in my different
approaches to multiple Dedekind zeta functions encouraged me to
continue my search of the right ones.  From Mladen Dimitrov I
learned about the cones constructed by Shintani. Anton Deitmar asked  many questions on the subject
that helped me to clear up (at least to myself) the structure of
the paper. Matthew Kerr clarified many questions I had about Hodge
theory. Dev Sinha explained simplicial de Rham structures on mapping spaces of a manifold, in particular, a path space or a loop space of a variety.  Finally, I would like to thank the referee whose effort greatly improved the presentation.

I acknowledge the great hospitality of
Durham University and the generous financial support from the
European Network Marie Currie   that allowed me to invite
Ronald Brown for several days.  I would like to thank Tubingen University and in
particular, Anton Deitmar and his grant ``Higher modular forms and
higher invariants" for the financial support and great working
conditions.


\renewcommand{\em}{\textrm}

\begin{small}

\renewcommand{\refname}{ {\flushleft\normalsize\bf{References}} }
    
\end{small}

\end{document}